\documentclass[12pt]{amsart}
\font\bbbld=msbm10 scaled\magstephalf

\newcommand{\bb}{\bar{b}}
\newcommand{\bi}{\bar{i}}
\newcommand{\bj}{\bar{j}}
\newcommand{\bk}{\bar{k}}
\newcommand{\bl}{\bar{l}}
\newcommand{\bm}{\bar{m}}
\newcommand{\bn}{\bar{n}}
\newcommand{\bp}{\bar{p}}
\newcommand{\bq}{\bar{q}}
\newcommand{\br}{\bar{r}}
\newcommand{\bs}{\bar{s}}
\newcommand{\bw}{\bar{w}}
\newcommand{\bz}{\bar{z}}

\newcommand{\bM}{\bar{M}}

\newcommand{\bpartial}{\bar{\partial}}

\def \p{\partial}
\def \f{\frac}

\def \o{\omega}
\def \d{\delta}

\def \m{\mu}

\def \tg{\tilde{g}}

\newcommand{\fg}{\mathfrak{g}}
\newcommand{\fRe}{\mathfrak{Re}}

\newcommand{\bfC}{\hbox{\bbbld C}}

\newcommand{\bfR}{\hbox{\bbbld R}}
\newcommand{\bfS}{\hbox{\bbbld S}}

\newcommand{\cA}{\mathcal{A}}

\newcommand{\cH}{\mathcal{H}}

\newcommand{\cT}{\mathcal{T}}

\newcommand{\ol}{\overline}
\newcommand{\ul}{\underline}

\newtheorem{theorem}{Theorem}[section]
\newtheorem{lemma}[theorem]{Lemma}
\newtheorem{proposition}[theorem]{Proposition}

\newtheorem{conjecture}[theorem]{Conjecture}
 \theoremstyle{definition}

\theoremstyle{remark}
\newtheorem{remark}[theorem]{Remark}

\numberwithin{equation}{section}

\begin{document}
\setlength{\baselineskip}{1.2\baselineskip}

\title[Complex Monge-Amp\`ere Equations]
{Complex Monge-Amp\`ere Equations \\ on Hermitian Manifolds}
\author{Bo Guan}
\address{Department of Mathematics, Ohio State University,
         Columbus, OH 43210}
\email{guan@math.osu.edu}
\thanks{Research of the first author was supported in part by NSF grants.}
\author{Qun Li}
\address{Department of Mathematics, Ohio State University,
         Columbus, OH 43210}
\email{li@math.osu.edu}

\begin{abstract}
We study complex Monge-Amp\`ere equations in Hermitian manifolds,
 both for the Dirichlet problem and in the case of compact manifolds
without boundary.
Our main results extend classical theorems of Yau~\cite{Yau78} and
Aubin~\cite{Aubin78} in the K\"ahler case, and those of
Caffarelli, Kohn, Nirenberg and Spruck~\cite{CKNS} for the Dirichlet
problem in $\bfC^n$.
As an application we study
the problem of finding geodesics in the space of Hermitian metrics,
generalizing existing results on Donaldson's conjecture \cite{Donaldson99}
in the  K\"ahler case.

\end{abstract}

\maketitle

\bigskip

\section{Introduction}
\label{gblq-I}
\setcounter{equation}{0}
\medskip

The complex Monge-Amp\`ere equation has close connections with
many important problems in complex geometry and analysis.
In the framework of K\"ahler geometry, it goes back at least to
Calabi~\cite{Calabi56} who conjectured that
any element in the first Chern class of a compact K\"ahler
manifold is the Ricci form of a K\"ahler metric cohomologous
to the underlying metric.
In \cite{Donaldson99}, Donaldson
made several conjectures on the space of K\"ahler metrics
which reduce to questions on a special Dirichlet problem for the
homogeneous complex Monge-Amp\`ere ({\em HCMA}) equation;
see also the related work of Mabuchi~\cite{Mabuchi87} and
Semmes~\cite{Semmes92}.
The HCMA equation, which is well defined on general complex manifolds,
also arises naturally in other interesting geometric problems.
One such example is the work on intrinsic norms by Chern,
Levine and Nirenberg~ \cite{CLN69}, Bedford and Taylor~\cite{BT79}
and P.-F. Guan~\cite{GuanPF02}, \cite{GuanPF08}.
There are also interesting results in the literature that connect the HCMA
equation on general complex manifolds with totally real submanifolds;
see, e.g. \cite{Wong}, \cite{Patrizio-Wong}, \cite{Lempert-Szoke},
\cite{Guillemin-Stenzel91} and \cite{Guillemin-Stenzel92}.

In ~\cite{Yau78}, Yau proved fundamental existence theorems of
classical solutions for complex Monge-Amp\`ere equations on compact K\"ahler
manifolds and consequently solved the Calabi conjecture. Yau's work
also shows the existence of K\"ahler-Einstein metrics on K\"ahler
manifolds with the first Chern number $c_1 (M) \leq 0$, confirming another
conjecture of Calabi~\cite{Calabi56} which was independently proved
by Aubin~\cite{Aubin78} in the case $c_1 (M) < 0$.
The classical solvability of the Dirichlet problem was established by
Caffarelli, Kohn, Nirenberg and Spruck~\cite{CKNS} for strongly pseudoconvex
domains in $\bfC^n$.
Later on the first author~\cite{Guan98b} extended their
results to general domains under the assumption of existence of
subsolutions.

Our primary goal in this paper is to attempt to extend these results to more
general geometric settings. We shall consider complex Monge-Amp\`ere equations
on Hermitian manifolds. Besides the technical challenges in the analytic aspect
which we shall discuss in more details later,
our motivation originates from trying to understand the above mentioned
Donaldson conjectures when one considers Hermitian metrics, as well as other
geometric problems some of which we shall treat in a forthcoming paper \cite{GL2}.

Let us first consider the Dirichlet problem.
Let $(M^n, \omega)$ be a compact Hermitian manifold of dimension
$n \geq 2$ with smooth boundary $\partial M$,  and $\bM = M \cup \partial M$.
Given $\psi \in C^{\infty} (M \times \bfR)$, $\psi > 0$,
and $\varphi \in C^{\infty} (\partial M)$, we seek solutions
of the complex Monge-Amp\`ere equation
\begin{equation}
\label{gblq-I10}
   \Big(\omega + \frac{\sqrt{-1}}{2} \partial \bpartial u\Big)^n
  = \psi (z, u) \omega^n \;\; \mbox{in $\bM$}
\end{equation}
satisfying the Dirichlet condition
\begin{equation}
\label{gblq-I20}
  u = \varphi \;\; \mbox{on $\partial M$}.
 \end{equation}

We require 
\begin{equation}
\label{gblq-I30}
\omega_u :=\omega + \frac{\sqrt{-1}}{2} \partial \bpartial u > 0
 \end{equation}
so that equation~\eqref{gblq-I10} is elliptic; we call such functions
{\em admissible}. Set
\begin{equation}
\label{gblq-I40}
 \cH = \{\phi \in C^2 (\bM): \omega_{\phi} > 0\}.
\end{equation}
We shall also call $\cH$ the {\em space of Hermitian metrics}.
As in the K\"ahler case, for $u \in \cH$,
$\omega_u$ is a Hermitian form on $M$ and equation~\eqref{gblq-I10}
describes one of its Ricci forms.

From the theory of fully nonlinear elliptic equations, a crucial step
in solving equation~\eqref{gblq-I10} is to derive {\em a priori} $C^2$
estimates for admissible solutions. Our first result is an extension of
Yau's estimate for $\Delta u$ \cite{Yau78} and the gradient estimates due
to Blocki~\cite{Blocki09} and P.-F. Guan~\cite{GuanPF} in the K\"ahler case.

\begin{theorem}
\label{gblq-th10}
Let $u \in \cH \cap C^4 (M)$ be a solution of equation~\eqref{gblq-I10}.
Then there exist positive constants $C_1$, $C_2$ depending on $|u|_{C^0 (\bM)}$
such that
\begin{equation}
\label{gblq-I50}
\max_{\bM} |\nabla u| \leq C_1 \Big(1 + \max_{\partial M} |\nabla u|\Big)
\end{equation}
and
\begin{equation}
\label{gblq-I60}
\max_{\bM} \Delta u\leq C_2 \Big(1 + \max_{\partial M} \Delta u\Big).
\end{equation}
\end{theorem}

More details are given in Propositions~\ref{gblq-prop-G10} and
\ref{gblq-prop-C10} of the dependence of $C_1$ and $C_2$ on $\psi$
and geometric quantities $M$ (torsion and curvatures).
Here we only emphasize that these constants do not depend on
$\inf \psi$ so the estimates \eqref{gblq-I50} and \eqref{gblq-I60} apply
to the degenerate case ($\psi \geq 0$).

The gradient estimate~\eqref{gblq-I50} was also proved independently by
Xiangwen Zhang~\cite{Zhang} who considered more general equations on compact
Hermitian manifolds without boundary.

A substantial difficulty in proving \eqref{gblq-I60} is to control
the extra terms involving third order derivatives which appear due to
the nontrivial torsion. To overcome this we construct a special local
coordinate system (Lemma~\ref{non-sym}) and make use of some unique properties
of the Monge-Amp\`ere operator. See the proof in Section~\ref{gblq-C2G}
for details.

In order to solve the Dirichlet problem~\eqref{gblq-I10}-\eqref{gblq-I20}
we also need estimates for second derivatives on the boundary.
For this we shall follow techniques developed in \cite{GS93},
\cite{Guan98a}, \cite{Guan98b} using subsolutions.
Our main existence result for the Dirichlet problem may be stated
as follows.

\begin{theorem}
\label{gblq-th20}
Suppose that there exists $\ul{u} \in C^0 (\bM)$ with $\omega_{\ul{u}} \geq 0$
 in $\bM$
(in weak sense \cite{BT76}), $\ul{u} = \varphi$ on $\partial M$
 and
\begin{equation}
\label{gblq-I10'}
  (\omega_{\ul{u}})^n \geq  \psi (z, \ul{u}) \omega^n \;\; \mbox{in $\bM$}.
 \end{equation}
Assume further that $\ul{u} \in C^2$ in a neighborhood of $\partial M$
(including $\partial M$). Then the Dirichlet
problem~\eqref{gblq-I10}-\eqref{gblq-I20} admits a solution
$u \in \cH \cap C^{\infty} (\bM)$ with $u \geq \ul{u}$ in $\bM$.
\end{theorem}

We note that Theorem~\ref{gblq-th10} also applies to compact manifolds
without boundary. Deriving the $C^0$ estimates, however, seems
a difficult question for general $\omega$. In the K\"ahler case,
Yau~\cite{Yau78} introduced a Moser iteration approach using his $C^2$
estimate and the Sobolev inequality. His proof was subsequently
simplified by Kazdan~\cite{Kazdan78} for $n =2$, and by Aubin~\cite{Aubin78}
and Bourguignon independently for arbitrary dimension (see e.g. \cite{Siu87}
and \cite{Tian00}).
Alternative proofs were given by Kolodziej~\cite{Kolodziej98} and
Blocki~\cite{Blocki05} based on the pluripotential theory (\cite{BT82}) and
the $L^2$ stability of the complex Monge-Amp\`ere operator (\cite{CP92}).
All these proofs seem to heavily rely on the closedness or, equivalently,
existence of local potentials of $\omega$ and it is not clear to us whether
any of them can be extended to the Hermitian case.
In this paper we impose the following condition
\begin{equation}
\label{gblq-I70}
\partial \bpartial \omega^k = 0, \;\; 1 \leq k \leq n-1
\end{equation}
which will also enable us to carry out the continuity method as in
\cite{Yau78}.

\begin{theorem}
\label{gblq-th30}
 Let $(M, \omega)$ be a compact Hermitian manifold without boundary
with $\omega$ satisfying \eqref{gblq-I70}.
Assume $\psi_u \geq 0$ and that there exists a function
$\phi \in C^{\infty} (M)$ such that
\begin{equation}
\label{gblq-I80}
\int_M \psi (z, \phi (z)) \omega^n  = \mbox{Vol}\, (M).
\end{equation}
Then there exists a solution $u \in \cH \cap C^{\infty} (M)$ of
equation~\eqref{gblq-I10}. Moreover the solution is unique, possibly up
to a constant.
\end{theorem}

Under stronger assumptions on $\psi$ condition \eqref{gblq-I70} may be
removed.

\begin{theorem}
\label{gblq-th40}
 Let $(M, \omega)$ be a compact Hermitian manifold without boundary.
If
\begin{equation}
\label{gblq-I90}
\lim_{u \rightarrow -\infty} \psi (\cdot, u) = 0,  \;\;
 \lim_{u \rightarrow +\infty} \psi (\cdot, u) = \infty.
\end{equation}
and $\psi_u > 0$, then equation~\eqref{gblq-I10} has a unique solution
in $\cH \cap C^{\infty} (M)$.
\end{theorem}

For applications in complex geometry it is very important to study the
degenerate complex Monge-Amp\`ere equation ($\psi \geq 0$ in \eqref{gblq-I10}). 
In general, the optimal regularity in the degenerate case is $C^{1,1}$;
see e.g., \cite{BF79}, \cite{GS80}, and there are many challenging 
open questions. In the forthcoming article~\cite{GL2} we shall focus on the
degenerate, especially the homogeneous, Monge-Amp\`ere equation in 
Hermitian manifolds and applications in geometric problems. In the current 
paper we shall only prove the following theorem for a special Dirichlet problem.

\begin{theorem}
\label{gblq-th50}
Let $M = N \times S$ where $N$ is a compact Hermitian manifold without
boundary, $\mbox{dim}_C N = n-1$,
and $S$ is a compact Riemann surface with smooth boundary $\partial S \neq \emptyset$.
Let $\omega$ be the product K\"ahler form on $M$.

Suppose that 
$\psi \geq 0$, $\psi^{\frac{1}{n}} \in  C^2 (\bM \times R)$, and
$\phi \in \cH \cap C^4 (\bM)$ satisfying
$(\omega_{\phi})^n \geq \psi$ on $\bM$.
Then there exists a weak admissible solution $u \in C^{1,\alpha} (\bM)$,
for all $\alpha \in (0, 1)$ with $\Delta u \in L^{\infty} (M)$,
of the Dirichlet problem 
  \begin{equation}
\label{gblq-I100}
 \left\{ \begin{aligned}
 & (\omega_u)^n = \psi
 \;\; \mbox{in $\bM$}, \\
 & u = \phi \;\; \mbox{on $\partial M$}. 
  \end{aligned} \right.
\end{equation}

Moreover, the solution is unique if $\psi_u \geq 0$, and $u \in C^{1,1} (\bM)$
provided that $M$ has nonnegative bisectional curvature.
\end{theorem}

As an immediate application of Theorem~\ref{gblq-th50} we can extend
existing results on a conjecture of Donaldson~\cite{Donaldson99} 
concerning geodesics in the space of K\"ahler metrics to the Hermitian 
setting; see Section~\ref{gblq-S} for details.

The paper is organized as follows. In Section~\ref{gblq-P} we
briefly recall some basic facts and formulas for Hermitian manifolds
and the complex Monge-Amp\`ere operator, fixing the notation along
the way. We shall also prove in this section the existence of local
coordinates with some special properties; see Lemma~\ref{non-sym}.
Such local coordinates are crucial to our proof of \eqref{gblq-I60}
in Section~\ref{gblq-C2G}, while the gradient estimate \eqref{gblq-I50}
is derived in Sections~\ref{gblq-G}.
Section~\ref{gblq-B} concerns the boundary estimates for second
derivatives. In Section~\ref{gblq-R} we come back to finish the
global estimates for all (real) second derivatives which enable us to
apply the Evans-Krylov theorem~\cite{Evans}, \cite{Krylov} to
obtain $C^{2, \alpha}$ and therefore higher order estimates by the classical
Schauder theory.
In Section~\ref{gblq-E} we discuss the $C^0$ estimates and
existence of solutions, completing the proof of
Theorems~\ref{gblq-th20}-\ref{gblq-th50}.
Finally, in Section~\ref{gblq-S} we extend results on the Donaldson conjecture
in the K\"ahler case to the space of Hermitian metrics.

The authors wish to thank Pengfei Guan and Fangyang Zheng for very helpful
conversations and suggestions.

\bigskip

\section{Preliminaries}
\label{gblq-P}
 \setcounter{equation}{0}
\medskip


 Let $M^n$ be a complex manifold of dimension $n$
and $g$ a Riemannian metric on $M$. Let $J$ be the induced
almost complex structure on $M$ so $J$ is integrable and
$J^2 = - \mbox{id}$. We assume that $J$ is compatible with $g$, i.e.
\begin{equation}
\label{cma-K10}
 g (u, v) = g (Ju, Jv), \;\; u, v \in TM;
\end{equation}
such $g$ is called a {\em Hermitian} metric.
Let $\omega$ be the {\em K\"ahler} form of $g$ defined
by
\begin{equation}
\label{cma-K20}
\omega (u, v) = - g (u, Jv).
\end{equation}
We recall that $g$ is {\em K\"ahler} if its K\"ahler form
$\omega$ is closed,  i.e. $d \omega = 0$.

The complex tangent bundle $T_C M = TM \times \bfC$ has a
natural splitting
\begin{equation}
\label{cma-K30}
 T_C M = T^{1,0} M + T^{0,1} M
\end{equation}
where $T^{1,0} M$ and  $T^{0,1} M$ are the $\pm
\sqrt{-1}$-eigenspaces of $J$. The metric $g$ is obviously extended
$\bfC$-linearly to $T_C M$, and
\begin{equation}
\label{cma-K40}
 g (u, v) = 0 \;\; \mbox{if $u,v \in T^{1,0} M$, or $u,v \in T^{0,1} M$}.
\end{equation}

Let $\nabla$ be the Chern connection of $g$. It satisfies
\begin{equation}
\label{cma-K55}
\nabla_u (g (v, w)) = g (\nabla_u v, w) + g (v, \nabla_u w)
\end{equation}
but may have nontrivial torsion.
The torsion $T$ and curvature $R$ of $\nabla$ are defined by
\begin{equation}
\label{cma-K95'}
T (u, v) = \nabla_u v - \nabla_v u - [u,v],
\end{equation}
and
\begin{equation}
\label{cma-K95}
 R (u, v) w = \nabla_u \nabla_v w - \nabla_v \nabla_u w - \nabla_{[u,v]} w,
 \end{equation}
 respectively.
Since $\nabla J = 0$ we have
\[ T (Ju, Jv) = J T (u, v) \]
and
\[ R (u, v) Jw = J R (u, v) w. \]
It follows that
\begin{equation}
\label{cma-K100'}
T (u, v, w) \equiv g (T (u, v), w) = g (T (Ju, Jv), Jw)
            \equiv T (Ju, Jv, Jw)
\end{equation}
and
\begin{equation}
\label{cma-K100}
R (u, v, w, x) \equiv g (R (u, v) w, x) = g (R (u, v) Jw, Jx).
\end{equation}
Therefore $R (u, v, w, x) = 0$ unless $w$ and $x$ are of different
type.

 In local coordinates $(z_1, \ldots, z_n)$, by \eqref{cma-K40}
\begin{equation}
\label{cma-K50}
 g \Big(\frac{\partial}{\partial z_j},
\frac{\partial}{\partial z_k}\Big) = 0, \;\;
  g \Big(\frac{\partial}{\partial \bz_j}, \frac{\partial}{\partial \bz_k}\Big)
   = 0
\end{equation}
since \
\begin{equation}
\label{cma-K60}
 J \frac{\partial}{\partial z_j} = \sqrt{-1}
\frac{\partial}{\partial z_j}, \;\;
 J \frac{\partial}{\partial \bz_j} = - \sqrt{-1}
\frac{\partial}{\partial \bz_j}.
\end{equation}
  We write
\begin{equation}
\label{cma-K70}
 g_{j\bk} = g \Big(\frac{\partial}{\partial z_j},
\frac{\partial}{\partial \bz_k}\Big),
  \;\; \{g^{i\bj}\} = \{g_{i\bj}\}^{-1}.
\end{equation}
That is, $g^{i\bj} g_{k\bj} = \delta_{ik}$.
 The K\"ahler form $\omega$ is then given by
\begin{equation}
\label{cma-K80}
 \omega = \frac{\sqrt{-1}}{2} g_{j\bk} dz_j \wedge d \bz_k.
\end{equation}

The Christoffel symbols $\Gamma^l_{jk}$ are defined by
\[ \nabla_{\frac{\partial}{\partial z_j}}
\frac{\partial}{\partial z_k}
   = \Gamma^l_{jk} \frac{\partial}{\partial z_l}. \]
Recall that by \eqref{cma-K55} and \eqref{cma-K50},
\begin{equation}
\label{cma-K85}
\left\{ \begin{aligned}
\nabla_{\frac{\partial}{\partial z_j}} \frac{\partial}{\partial \bz_k}
 \,& = \nabla_{\frac{\partial}{\partial \bz_j}} \frac{\partial}{\partial z_k}
     = 0, \\
 \nabla_{\frac{\partial}{\partial \bz_j}} \frac{\partial}{\partial \bz_k}
 \,& = \Gamma^{\bl}_{\bj \bk} \frac{\partial}{\partial \bz_l}
     = \ol{\Gamma^l_{jk}} \frac{\partial}{\partial \bz_l}
 \end{aligned} \right.
\end{equation}
and
\begin{equation}
\label{cma-K90}
\Gamma^l_{jk} = g^{l\bm} \frac{\partial g_{k\bm}}{\partial z_j}.
\end{equation}
For the torsion and curvature we use the standard notion
\[ T_{ijk} = T \Big(\frac{\partial}{\partial z_i},
   \frac{\partial}{\partial z_j}, \frac{\partial}{\partial z_k}\Big),
   \;\; \mbox{etc.} \]
and
\[ R_{i\bj k\bl} = R \Big(\frac{\partial}{\partial z_i},
 \frac{\partial}{\partial \bz_j}, \frac{\partial}{\partial z_k},
 \frac{\partial}{\partial \bz_l}\Big), \;\; \mbox{etc.} \]
Obviously,
\[ T_{ijk} = T_{i\bj k}  = T_{i\bj \bk} = 0, \\
   T_{ij\bk} = \frac{\partial g_{j\bk}}{\partial z_i}
               - \frac{\partial g_{i\bk}}{\partial z_j} \]
and
\begin{equation}
\label{cma-K102}
 T^k_{ij} \equiv g^{k\bl}  T_{ij\bl}
   = \Gamma^k_{ij} - \Gamma^k_{ji} = g^{k\bl} \Big(\frac{\partial
g_{j\bl}}{\partial z_i} - \frac{\partial g_{i\bl}}{\partial z_j}\Big).
\end{equation}
From \eqref{cma-K100} and \eqref{cma-K90} it follows that
\begin{equation}
\label{cma-K105}
R_{i\bj k l}  = R_{i j k l} = 0,
\end{equation}
\begin{equation}
\label{cma-K110}
\begin{aligned}
 R_{i\bj k\bl} 
 = - g_{m \bl} \frac{\partial \Gamma_{ik}^m}{\partial \bz_j}
 = - \frac{\partial^2 g_{k\bl}}{\partial z_i \partial \bz_j}
       + g^{p\bq} \frac{\partial g_{k\bq}}{\partial z_i}
                  \frac{\partial g_{p\bl}}{\partial \bz_j}
         \end{aligned}
\end{equation}
and
\begin{equation}
\label{cma-K115'}
\begin{aligned}
 R_{i j k\bl} 
 = g_{m \bl} \Big(\frac{\partial \Gamma_{jk}^m}{\partial z_i}
                   - \frac{\partial \Gamma_{ik}^m}{\partial z_j}
   + \Gamma_{iq}^m \Gamma_{jk}^q - \Gamma_{jq}^m \Gamma_{ik}^q\Big).
\end{aligned}
\end{equation}
By \eqref{cma-K110} and \eqref{cma-K102} we have
\begin{equation}
\label{cma-K115}
R_{i \bj k\bl} - R_{k \bj i\bl}
   =  g_{m \bl} \frac{\partial T_{ki}^m}{\partial \bz_j}
   = g_{m \bl} \nabla_{\bj} T_{ki}^m,
\end{equation}
which also follows form the general Bianchi identity.

The traces of the curvature tensor
\begin{equation}
\label{cma-K120'}
 R_{k\bl} = g^{i\bj} R_{i\bj k\bl}
   = - g^{i\bj} \frac{\partial^2 g_{k\bl}}{\partial z_i \partial \bz_j}
       + g^{i\bj} g^{p\bq} \frac{\partial g_{k\bq}}{\partial z_i}
                  \frac{\partial g_{p\bl}}{\partial \bz_j}
\end{equation}
and
\begin{equation}
\label{cma-K120"}
 S_{i\bj} \equiv g^{k\bl} R_{i\bj k\bl}
     = - g^{k\bl} \frac{\partial^2 g_{k\bl}}{\partial z_i \partial \bz_j}
       + g^{k\bl} g^{p\bq} \frac{\partial g_{k\bq}}{\partial z_i}
                  \frac{\partial g_{p\bl}}{\partial \bz_j}
 \end{equation}
are called the {\em first} and {\em second} Ricci tensors, respectively.
Note that
\begin{equation}
\label{cma-K120}
 S_{i\bj} = -\frac{\partial^2}{\partial z_i \partial \bz_j} \log\det g_{k\bl}.
\end{equation}

The following special local coordinates will be crucial to our proof of
the \textit{a priori} estimates for $\Delta u$ in Theorem~\ref{gblq-th10}.

\begin{lemma}
\label{non-sym}
Around a point $p \in M$ there exist local coordinates such that, at $p$,
\begin{equation}
\label{coord}
\begin{aligned}
   & g_{i\bj} = \delta_{ij},   \;\;
    \frac{\partial g_{i\bi}}{\partial z_j} = 0, \;\; \forall \;i, j.
\end{aligned}
\end{equation}
\end{lemma}

\begin{proof}
Let $(z_1, z_2, \cdots, z_n)$ be a local coordinate system around $p$
such that $z_i(p)=0$ for $i=1, \cdots, n$ and
\[ g_{i\bj}(p) := g \Big(\f {\p}{\p z_i}, \f {\p}{\p \bar{z_j}}\Big)
    = \d_{i j}. \]
Define new coordinates $(w_1, w_2, \cdots, w_n)$ by
\begin{equation}
\label{new-coord}
w_r = z_r + \sum_{m \neq r} \f {\p g_{r\br}}{\p z_m}(p)z_m z_r
          + \frac{1}{2} \f {\p g_{r\br}}{\p z_r}(p) z_r^2,
\;\; 1 \leq r \leq n.
\end{equation}
We have
\begin{equation}\label{new-metric}
\tg_{i\bj}: =g \Big(\f {\p}{\p w_i}, \f {\p}{\p \bar{w_j}}\Big)
 =\sum_{r,s} g_{r\bs} \f {\p z_r}{\p w_i} \overline{\f {\p z_s}{\p w_j}}.
\end{equation}
It follows that
\begin{equation}
\label{ijj}
\f {\p \tg_{i\bj}}{\p w_k} =
  \sum_{r,s} g_{r\bs} \f {\p^2 z_r}{\p w_i \p w_k}
    \overline{\f {\p z_s}{\p w_j}}
  + \sum_{r,s, p} \f {\p g_{r \bs}}{\p z_p}
    \f {\p z_p}{\p w_k} \f {\p z_r}{\p w_i}  \overline{\f {\p z_s}{\p w_j}}.
\end{equation}

Differentiate (\ref{new-coord}) with respect to $w_i$ and $w_k$.
We see that, at $p$,
\[ \f {\p z_r}{\p w_i}=\d_{r i}, \;\;
\f {\p^2 z_r}{\p w_i \p w_k} = - \sum_{m \neq r} \f {\p g_{r\br}}{\p z_m}
   \Big(\f {\p z_m}{\p w_i} \f {\p z_r}{\p w_k}
      + \f {\p z_m}{\p w_k}\f {\p z_r}{\p w_i}\Big)
      - \f {\p g_{r\br}}{\p z_r} \f {\p z_r}{\p w_i} \f {\p z_r}{\p w_k}. \]
Plugging these into (\ref{ijj}), we obtain at $p$,
\begin{equation}
\label{gblq-P80}
\left\{
\begin{aligned}
\frac{\partial \tg_{i\bi}}{\partial w_k}
   = \, & \frac{\partial g_{i\bi}}{\partial z_k}
      - \frac{\partial g_{i\bi}}{\partial z_k} = 0, \;\; \forall \; i, k,\\
\frac{\partial \tg_{i\bj}}{\partial w_j}
   = \,& \frac{\partial g_{i\bj}} {\partial z_j}
      - \frac{\partial g_{j\bj}}{\partial z_i} = T_{ji}^j,
      \;\; \forall \; i \neq j, \\
\frac{\partial \tg_{i\bj}}{\partial w_k}
   = \,& \frac{\partial g_{i\bj} }{\partial z_k}, \;\; \mbox{otherwise.}
\end{aligned}
\right.
\end{equation}
Finally, switching $w$ and $z$ gives \eqref{coord}.
\end{proof}

\begin{remark}
If, in place of \eqref{new-coord}, we
define
\begin{equation}
\label{new-coord'}
w_r = z_r + \sum_{m \neq r} \f {\p g_{m\br}}{\p z_r}(p)z_m z_r
          + \frac{1}{2} \f {\p g_{r\br}}{\p z_r}(p) z_r^2,
\;\; 1 \leq r \leq n,
\end{equation}
then under the new coordinates $(w_1, w_2, \cdots, w_n)$,
\begin{equation}
\label{gblq-P80'}
\left\{
\begin{aligned}
\frac{\partial \tg_{i\bj}}{\partial w_j} (p)
   = \, & 0, \;\; \forall \; i, j,\\
\frac{\partial \tg_{i\bi}}{\partial w_k} (p)
      = \,& T_{ki}^i, \;\; \forall \; i \neq k, \\
\frac{\partial \tg_{i\bj}}{\partial w_k} (p)
   = \,& \frac{\partial g_{i\bj} }{\partial z_k} (p), \;\; \mbox{otherwise.}
 \end{aligned}
\right.
\end{equation}
\end{remark}

The following lemma and its proof can be found in \cite{ST}.

\begin{lemma}
\label{lemma-I10}
Around a point $p \in M$ there exist local coordinates such that, at $p$,
\begin{equation}
\label{gblq-G100}
\left\{
\begin{aligned}
   & g_{i\bj} = \delta_{ij},   \\
   & \frac{\partial g_{i\bj}}{\partial z_k}
      + \frac{\partial g_{k\bj}}{\partial z_i} = 0.
\end{aligned}
\right.
\end{equation}
Consequently, $T_{ij}^k = 2 \frac{\partial g_{j\bk}}{\partial z_i}$ at $p$.
\end{lemma}

\begin{remark}
In general it is impossible to find local coordinates satisfying
both \eqref{coord} and \eqref{gblq-G100} simultaneously.
\end{remark}

Let $\Lambda^{p, q}$ denote differential forms of type $(p, q)$ on $M$.
The exterior differential $d$ has a natural decomposition
$d = \partial + \bpartial$ where
\[ \partial: \Lambda^{p, q} \rightarrow \Lambda^{p+1, q} , \;\;
  \bpartial: \Lambda^{p, q} \rightarrow \Lambda^{p, q+1}.  \]
Recall that $\partial^2 = \bpartial^2 = \partial \bpartial + \bpartial \partial = 0$
and, by the Stokes theorem
\[ \int_M \partial \alpha = \int_{\partial M} \alpha,
 \;\; \forall \; \alpha \in \Lambda^{n-1, n}. \]
A similar formula holds for $\bpartial$.

For a function $u \in C^2 (M)$, $\partial \bpartial u$
is given in local coordinates by
\begin{equation}
\label{cma-K80'}
\partial \bpartial u = \frac{\partial^2 u}{\partial z_i \partial \bz_j}
  dz_i \wedge d \bz_j.
\end{equation}
We use $\nabla^2 u$ to denote the {\em Hessian} of $u$:
\begin{equation}
\label{cma-K215}
\nabla^2 u (X, Y) \equiv \nabla_Y \nabla_X u
   = Y (X u) - (\nabla_Y X) u, \;\; X, Y \in TM.
\end{equation}
By \eqref{cma-K85} we see that
\begin{equation}
\label{cma-K220}
\nabla_{\frac{\partial}{\partial z_i}}
 \nabla_{\frac{\partial}{\partial \bz_j}} u
   = \frac{\partial^2 u}{\partial z_i \partial \bz_j}.
\end{equation}
Consequently, the Laplacian of $u$ with respect to the Chern connection is
\begin{equation}
\label{cma-K230}
\Delta u = g^{i\bj} \frac{\partial^2 u}{\partial z_i \partial \bz_j},
\end{equation}
or equivalently,
\begin{equation}
\label{cma-K235}
\Delta u \omega^n = \frac{\sqrt{-1}}{2} \partial \bpartial u \wedge \omega^{n-1}.
\end{equation}
Integrating \eqref{cma-K235} (by parts), we obtain
\begin{equation}
\label{cma-K238}
\begin{aligned}
\frac{2}{\sqrt{-1}} \int_M \Delta u \omega^n
  = \,& \int_M \partial \bpartial u \wedge \omega^{n-1} \\
  = \,& \int_{\partial M} \bpartial u \wedge \omega^{n-1}
        + \int_M  \bpartial u \wedge \partial \omega^{n-1} \\
  = \,& \int_{\partial M} (\bpartial u \wedge \omega^{n-1} + u \partial \omega^{n-1})
   + \int_M  u \partial \bpartial \omega^{n-1}.
   \end{aligned}
\end{equation}

Finally, in the following sections where we derive the {\em a priori} estimates
we shall consider the slightly more general equation
\begin{equation}
\label{gblq-I10mu}
 \det \Big(\mu \omega + \frac{\sqrt{-1}}{2} \partial \bpartial u\Big)^n
    = \psi (z, u) \omega^n
\end{equation}
where $\mu$ is a given smooth function on $M$ and may also depend on $u$ and $\nabla u$.
We shall abuse the notation to write
\[ \omega_u := \mu \omega + \frac{\sqrt{-1}}{2} \partial \bpartial u
       \in \Lambda^{1, 1}, \]
and $u \in \cH$ means $\omega_u > 0$. In local coordinates
 equation~\eqref{gblq-I10mu} takes the form
\begin{equation}
\label{cma2-M10}
 \det \Big(\mu g_{i\bj} + \frac{\partial^2 u}{\partial z_i \partial \bz_j}\Big)
    = \psi (z, u) \det g_{i\bj}.
\end{equation}
In the current paper we assume $\mu= \mu (z, u)$
 and $|\mu| > 0$ in $\bM \times \bfR$. In \cite{GL2}
 we shall consider other cases including $\mu = 0$.

\bigskip

\section{Gradient estimates}
\label{gblq-G}
\setcounter{equation}{0}
\medskip

In this section we derive the gradient estimate \eqref{gblq-I50}
for a solution $u \in \cH \cap C^3 (M)$ of (\ref{gblq-I10mu}).
 Throughout this and next sections we shall use ordinary derivatives.
For convenience we write in local coordinates,
\[ u_i = \frac{\partial u}{\partial z_i}, \;\;
   u_{\bi} = \frac{\partial u}{\partial \bz_i}, \;\;
   u_{i\bj} = \frac{\partial^2 u}{\partial z_i \partial \bz_j}, \;\;
   g_{i\bj k} = \frac{\partial g_{i\bj}}{\partial z_k},  \;\;
   g_{i\bj k\bl} = \frac{\partial^2 g_{i\bj}}{\partial z_k \partial \bz_l},  \;\;
\mbox{etc}, \]
and
\[ \fg_{i\bj} = u_{i\bj} + \mu g_{ij}, \;\; \{\fg_{i\bj}\} = \{\fg_{i\bj}\}^{-1}. \]

We first present some calculation. In local coordinates,
\begin{equation}
\label{gblq-G10}
|\nabla u|^2 = g^{k\bl} u_k u_{\bl}.
\end{equation}
We have
\begin{equation}
\label{cma2-M50'}
 (|\nabla u|^2)_i = g^{k\bl} (u_{ki} u_{\bl} + u_k u_{i\bl})
         - g^{k\bb} g^{a\bl} g_{a\bb i} u_k u_{\bl}
\end{equation}
\begin{equation}
\label{cma2-M60'}
\begin{aligned}
(|\nabla u|^2)_{i\bj} = \,
  & g^{k\bl} (u_{k \bj} u_{i\bl} + u_{i\bj k} u_{\bl}
     + u_k u_{i \bj \bl} + u_{ki} u_{\bl \bj}) \\
  &  - g^{k\bb} g^{a\bl} g_{a\bb i} (u_{k \bj} u_{\bl} + u_k u_{\bl \bj})
     - g^{k\bb} g^{a\bl} g_{a\bb \bj} (u_{ki} u_{\bl} + u_k u_{i\bl}) \\
  & + [(g^{k\bq} g^{p\bb} g^{a\bl} + g^{k\bb} g^{a\bq} g^{p\bl}) g_{a\bb i}
      g_{p\bq \bj} - g^{k\bb} g^{a\bl} g_{a \bb i \bj}] u_k u_{\bl}.
\end{aligned}
\end{equation}
Here we have used the formula
\[ (g^{k\bl})_i = - g^{k\bb} g^{a\bl} g_{a\bb i}. \]

Let $p$ be a fixed point. We may assume that
\eqref{gblq-G100} holds at $p$ and
\begin{equation}
\label{gblq-G110}
 \mbox{$\{u_{i\bj}\}$ is diagonal.}
\end{equation}
 Thus
\begin{equation}
\label{cma2-M50}
 (|\nabla u|^2)_i =  u_i u_{i\bi} + (u_{ki}  - g_{k\bl i} u_{l}) u_{\bk}
\end{equation}
and
\begin{equation}
\label{cma2-M60}
\begin{aligned}
(|\nabla u|^2)_{i\bi} = \,
  & u_{i\bi}^2 + u_{i\bi k} u_{\bk} + u_{i \bi \bk} u_k
     + \sum_k |u_{ki} - g_{k\bl i} u_{l}|^2  \\
  &  - 2 \fRe \{g_{l\bi i} u_{i \bi} u_{\bl}\}
     + (g_{l\bp i} g_{p\bk \bi} - g_{l \bk i \bi}) u_k u_{\bl}.
\end{aligned}
\end{equation}

\begin{lemma}
\label{lemma-PL20'}
Let $\phi$ be a function such that $e^{\phi} |\nabla u|^2$ attains its
maximum at an interior point $p$ where, in local coordinates,
\eqref{gblq-G100} and \eqref{gblq-G110} hold. Then, at $p$,
\begin{equation}
\label{bglq-M20}
\begin{aligned}
 |\nabla u|^2 \sum \fg^{i\bi} \phi_{i\bi} \leq
 \, &  - 2  \sum \fRe \{\fg^{i\bi} u_i u_{i\bi} \phi_{\bi}\}
       + 2 |\nabla f| |\nabla u| - 2 |\nabla u|^2 f_u \\
    &  - |\nabla u|^2 \Big(\inf_{j, l} R_{j\bj l\bl}
       - \sup_{j, l} |T_{lj}^j|^2\Big) \sum \fg^{i\bi} \\
    &  + 2 \fRe \{(\mu)_{k}  u_{\bk}\} \sum \fg^{i\bi}.
  \end{aligned}
\end{equation}
\end{lemma}

\begin{proof}
In the proof all calculations are done at $p$ where
since $e^{\phi} |\nabla u|^2$ attains its maximum,
\begin{equation}
\label{gblq-G30}
  \frac{(|\nabla u|^2)_i}{|\nabla u|^2} + \phi_i = 0, \;\;
\frac{(|\nabla u|^2)_{\bi}}{|\nabla u|^2} + \phi_{\bi} = 0
\end{equation}
and
\begin{equation}
\label{gblq-G40}
 \frac{(|\nabla u|^2)_{i\bi}}{|\nabla u|^2}
-  \frac{|(|\nabla u|^2)_{i}|^2}{|\nabla u|^4}
  + \phi_{i\bi} \leq 0.
\end{equation}
By \eqref{cma2-M50} and \eqref{gblq-G30},
\begin{equation}
\label{gblq-G50}
 |(|\nabla u|^2)_i|^2 = \sum_k |u_{k}|^2 |u_{ki}  - g_{k\bl i} u_{l}|^2
 - 2 |\nabla u|^2 \fRe \{u_i u_{i\bi} \phi_{\bi}\} - |u_i|^2  u_{i\bi}^2.
\end{equation}

Differentiating equation (\ref{cma2-M10}) we have
\begin{equation}
\label{gblq-M70}
\begin{aligned}
\fg^{i\bj} (u_{i\bj k} u_{\bk} + u_k u_{i\bj \bk})
    = & \, 2 |\nabla u|^2 f_u + 2 \fRe \{f_{z_k} u_{\bk} + g^{i\bj} g_{i\bj k} u_{\bk}\} \\
      & - 2 \fRe \{\fg^{i\bj} (\mu g_{i\bj k} + (\mu)_{k} g_{i\bj}) u_{\bk}\}.
       \end{aligned}
\end{equation}
Note that by \eqref{gblq-G100},
\begin{equation}
\label{gblq-M75}
\sum_{i,l} g_{i \bi l}  u_{\bl} - \mu \fg^{i\bi} g_{i \bi l} u_{\bl}
 = \fg^{i\bi} u_{i \bi} g_{i \bi l} u_{\bl}
 = - \fg^{i\bi} u_{i \bi} g_{l \bi i} u_{\bl}.
\end{equation}
From (\ref{cma2-M60}) and (\ref{gblq-M70}) we see that
\begin{equation}
\label{bglq-M80}
\begin{aligned}
\fg^{i\bi} (|\nabla u|^2)_{i\bi} =
\, & \fg^{i\bi} |u_{i\bi} - 2  g_{l \bi i} u_{\bl}|^2
     + \sum_{i, k} \fg^{i\bi} |u_{ki} - g_{k\bl i} u_{l}|^2  \\
   & + \fg^{i\bi} (g_{l\bp i} g_{p\bk \bi} - g_{l \bk i \bi}) u_k u_{\bl}
     - 4 \fg^{i\bi} |g_{l \bi i} u_{\bl}|^2 \\
   & + 2 |\nabla u|^2 f_u + 2 \fRe \{f_{z_k} u_{\bk}\}
     - 2 \fRe \{(\mu)_{k}  u_{\bk}\} \sum \fg^{i\bi}.
\end{aligned}
\end{equation}
Combining \eqref{gblq-G40}, \eqref{gblq-G50} and \eqref{bglq-M80}, we obtain
\begin{equation}
\label{cma2-M40}
\begin{aligned}
0 \geq \,
&  |\nabla u|^2 \Big(\inf_{i,l} \{g_{l\bp i} g_{p\bl \bi} - g_{l \bl i \bi}\}
- 4 \sup_{i,l} |g_{l \bi i}|^2 \Big)
    \sum \fg^{i\bi} \\
&  + 2 |\nabla u|^2 f_u - 2 |\nabla f| |\nabla u|
   - 2 \fRe \{(\mu)_{k}  u_{\bk}\} \sum \fg^{i\bi} \\
&  + 2 \fRe \sum \{\fg^{i\bi} u_i u_{i\bi} \phi_{\bi}\}
   + |\nabla u|^2 \sum \fg^{i\bi} \phi_{i\bi}.
  \end{aligned}
\end{equation}
This proves \eqref{bglq-M20}.
\end{proof}

\begin{proposition}
\label{gblq-prop-G10}
There exists $C > 0$ depending on
\[ \sup_{M} |u|, \; \inf_{M} (\psi^{\frac{1}{n}})_u, \; \sup_{M} |\nabla \psi^{\frac{1}{n}}|, \;
 \inf_M \frac{1}{|\mu|} \Big(\inf_{j, l} R_{j\bj l\bl}
       - \sup_{j, l} |T_{lj}^j|^2\Big), \]
       and
 \[  \sup_{M} \{|\mu|^{-1} + |\nabla \log |\mu|| + (\log |\mu|)_u\} \]
such that
\begin{equation}
\label{cma-40}
\max_{\bM} |\nabla u| \leq 
 C (1 + \max_{\partial M}  |\nabla u|).  
\end{equation}
\end{proposition}

\begin{proof}
Let $L = \inf_M u$ and $\phi =A e^{\nu (L-u)}$ where $A > 0$ to be
determined later and $\nu = \mu/|\mu|$. We have
\begin{equation}
\label{gblq-G60}
\begin{aligned}
 e^{\nu (u-L)} \fRe\{\fg^{i\bi} u_i u_{i\bi} \phi_{\bi}\}
  =  - \nu A |\nabla u|^2 + |\mu| A \sum \fg^{i\bi} u_{i} u_{\bi}
     \end{aligned}
\end{equation}
and
\begin{equation}
\label{gblq-G70}
\begin{aligned}
 e^{\nu (u-L)} \sum \fg^{i\bi} \phi_{i\bi}
  = & \, A  \fg^{i\bi} (u_{i} u_{\bi} - \nu u_{i\bi}) \\
  = & \, A \fg^{i\bi} u_{i} u_{\bi} -  n \nu A
   + |\mu| A \sum \fg^{i\bi} \\
 \geq & \, A |\nabla u|^2 \min_{i} \fg^{i\bi} + |\mu| A \sum \fg^{i\bi}
          - n \nu A\\
\geq & \, n A (|\mu|^{\frac{n-1}{n}} \psi^{-\frac{1}{n}}
                 |\nabla u|^{\frac{2}{n}} - \nu).
     \end{aligned}
\end{equation}
 By Lemma~\ref{lemma-PL20'}, at an interior point where $e^{\phi} |\nabla u|^2$
 achieves its maximum we have
\begin{equation}
\label{bglq-M20'}
\begin{aligned}
A e^{\nu (L-u)}  |\nabla u|^2 & \, \Big(|\mu|^{\frac{n-1}{n}} \psi^{-\frac{1}{n}}
                 |\nabla u|^{\frac{2}{n}} - 3 \nu
+ \frac{(n-1)|\mu|}{n} \sum \fg^{i\bi}\Big) \\
\leq \,
    &  2 |\nabla f| |\nabla u| - 2 |\nabla u|^2 f_u
       + 2 (\mu_u |\nabla u|^2 + |\nabla \mu| |\nabla u|) \sum \fg^{i\bi} \\
    &  - |\nabla u|^2 \Big(\inf_{j, l} R_{j\bj l\bl}
       - \sup_{j, l} |T_{lj}^j|^2\Big) \sum \fg^{i\bi}.
\end{aligned}
\end{equation}
Choose $A$ such that
\[ \frac{n-1}{n} |A| \geq \sup_M \frac{e^{u-L}} {|\mu|} \Big(1 + 2 \mu_u
    - \inf_{j, l} R_{j\bj l\bl}
       + \sup_{j, l} |T_{lj}^j|^2\Big). \]
We obtain
\[ |\nabla u| \leq 2 |\nabla \log |\mu|| e^{u-L} \]
or
\[ A |\mu|^{\frac{n-1}{n}} |\nabla u|^{\frac{n+2}{n}}
   - 3 \psi^{\frac{1}{n}} |\nabla u|
  + 2 n e^{u-L} ((\psi^{\frac{1}{n}})_u |\nabla u|
   - |\nabla \psi^{\frac{1}{n}}|) \leq 0.
 \]
 This gives the estimate in \eqref{cma-40}.
\end{proof}

\bigskip

\section{Global estimates for $\Delta u$}
\label{gblq-C2G}
\setcounter{equation}{0}
\medskip

 In this section we derive the estimate \eqref{gblq-I60}.
We wish to include the degenerate case $\psi \geq 0$.
So we shall still assume $\psi > 0$ but
the estimates will not depend on the lower bound of $\psi$.
We shall follow the notations in Section~\ref{gblq-G} and use
ordinary derivatives.

Throughout Sections~\ref{gblq-C2G}-\ref{gblq-R} we assume
$u$  is a solution of \eqref{gblq-I10} in $\cH \cap C^4 (M)$.
  In local coordinates,
\[ \Delta u  = g^{k\bl} u_{k\bl}.  \]
Therefore,
\begin{equation}
\label{gblq-C10}
 (\Delta u)_i = g^{k\bl} u_{k\bl i}
         - g^{k\bb} g^{a\bl} g_{a\bb i} u_{k \bl}
\end{equation}
\begin{equation}
\label{gblq-C20}
\begin{aligned}
(\Delta u)_{i\bi} = \,
  & g^{k\bl} u_{k\bl i \bi}
    - g^{k\bb} g^{a\bl} (g_{a\bb i} u_{k \bl \bi} + g_{a\bb \bi} u_{k \bl i}) \\
  & + [(g^{k\bq} g^{p\bb} g^{a\bl} + g^{k\bb} g^{a\bq} g^{p\bl}) g_{a\bb i}
      g_{p\bq \bi} - g^{k\bb} g^{a\bl} g_{a \bb i \bi}] u_{k\bl}.
\end{aligned}
\end{equation}

\begin{lemma}
\label{gblq-lemma-C10}
Assume that \eqref{coord} and \eqref{gblq-G110} hold at $p \in M$ .
Then at $p$,
\begin{equation}
\label{gblq-C120}
\begin{aligned}
\fg^{i\bi} (\Delta u)_{i\bi} \geq
\, &  \fg^{i\bi} \fg^{j\bj} |u_{i\bj j} + (\mu g_{i\bj})_j|^2
       + \Delta (f) - 2 \fg^{i\bi} \fRe \{T_{ik}^i (\mu)_{\bk}\} \\
   &   - (|\mu| |T|^2 + \Delta (\mu)) \sum \fg^{i\bi}
       - \mu \fg^{i\bi} (S_{i\bi} - R_{i\bi}) \\
   &   + \inf_{j,k} R_{j\bj k\bk} \Big((\Delta u + n \mu) \sum \fg^{i\bi} - n^2\Big).
       \end{aligned}
\end{equation}
\end{lemma}

\begin{proof}
By \eqref{gblq-C10} and \eqref{gblq-C20},
\begin{equation}
\label{gblq-C30}
 (\Delta u)_i =  u_{k\bk i} -  g_{k\bk i} u_{k \bk},
\end{equation}
\begin{equation}
\label{gblq-C40}
\begin{aligned}
(\Delta u)_{i\bi} = \,
  & u_{i\bi k\bk} - 2 \fRe \{u_{k \bj i} g_{j \bk \bi}\}
    + (g_{k\bp i} g_{p\bk \bi} + g_{p\bk i} g_{k\bp \bi}
        -  g_{k \bk i \bi}) u_{k\bk} \\
 = \, & u_{i\bi k\bk} - 2 \fRe \{u_{k \bj i} g_{j \bk \bi}\}
    + (\fg_{k\bk} - \mu) (g_{p\bk i} g_{k\bp \bi} + R_{i\bi k\bk}).
\end{aligned}
\end{equation}
Differentiating equation (\ref{cma2-M10}) twice we obtain
\begin{equation}
\label{gblq-C70}
\begin{aligned}
\fg^{i\bi} u_{i\bi k\bk}  = \,
        & \fg^{i\bi} \fg^{j\bj} |u_{i\bj k} + (\mu g_{i\bj})_k|^2
          + (f)_{k\bk} + g_{i\bi k \bk} - g_{i\bj k} g_{j\bi \bk}
          - \fg^{i\bi} (\mu g_{i\bi})_{k \bk} \\
  = \, &  \fg^{i\bi} \fg^{j\bj} |u_{i\bj k} + (\mu g_{i\bj})_k|^2
         + (f)_{k\bk} - R_{k\bk i\bi}
         -  (\mu)_{k \bk}  \sum \fg^{i\bi} \\
       & + \mu \fg^{i\bi} (R_{k\bk i\bi} - |g_{i\bj k}|^2)
         -  2 \fg^{i\bi} \fRe\{(\mu)_k g_{i\bi \bk}\}.
       \end{aligned}
\end{equation}
From \eqref{coord} we have,
\begin{equation}
\label{gblq-C100}
 \sum_{j, k} u_{k \bj i} g_{j \bk \bi}
 = \sum_{j \neq k} [u_{i \bj k} + (\mu g_{i\bj})_k] g_{j \bk \bi}
 - \mu \sum_{j \neq k} g_{i\bj k} g_{j \bk \bi} - \sum_k (\mu)_k g_{i\bk \bi}.
\end{equation}
By Cauchy-Schwarz inequality,
\begin{equation}
\label{gblq-C110}
\begin{aligned}
 2 \sum_{j \neq k} |\fRe \{[u_{i \bj k} + (\mu g_{i\bj})_k] g_{j \bk \bi}\}|
 \leq \, & \sum_{j \neq k} \fg^{j\bj} |u_{i \bj k} + (\mu g_{i\bj})_k|^2
 + \sum_{j \neq k} \fg_{j\bj} |g_{k\bj i}|^2.
   \end{aligned}
\end{equation}
Finally, combining  \eqref{gblq-C40},  \eqref{gblq-C70}, \eqref{gblq-C100}, 
\eqref{gblq-C110} and
\[ |g_{i\bj k}|^2 + |g_{k\bj i}|^2 - 2 \fRe \{g_{i\bj k} g_{j \bk \bi}\}
  = |g_{k\bj i} - g_{i\bj k}|^2 = |T_{ik}^j|^2 \]
we derive
\begin{equation}
\label{gblq-C120'}
\begin{aligned}
\fg^{i\bi} (\Delta u)_{i\bi} \geq
\, &  \fg^{i\bi} \fg^{j\bj} |u_{i\bj j} + (\mu g_{i\bj})_j|^2
     + \Delta (f) + (\fg^{i\bi} \fg_{j\bj} - 1) R_{i\bi j\bj} \\
   & - \mu \sum_k \fg^{i\bi} (R_{i\bi k\bk} - R_{k\bk i\bi})
     - \mu \sum_{j,k} \fg^{i\bi} |T_{ik}^j|^2 \\
   & - 2 \fg^{i\bi} \fRe \{T_{ik}^i (\mu)_{\bk}\}
     - \Delta (\mu) \sum \fg^{i\bi}.
       \end{aligned}
\end{equation}
which gives \eqref{gblq-C120}. 
\end{proof}

\begin{lemma}
\label{gblq-lemma-C20}
Suppose $e^\phi (n \mu + \Delta u)$ achieves its maximum at an interior point
$p \in M$ where \eqref{coord} and \eqref{gblq-G110} hold. Then, at $p$,
\begin{equation}
\label{gblq-C150}
  \begin{aligned}
 (n \mu + \Delta u) \fg^{i\bi} \phi_{i\bi}
      \,& +  2 \fg^{i\bi} \fRe\{\phi_{i} \bar{\lambda_{i}}\} \\
 \leq \,& - (n \mu + \Delta u) \inf_{j,k} R_{j\bj k\bk} \sum \fg^{i\bi}
           - \Delta (f) \\
        &  + n^2 \inf_{j,k} R_{j\bj k\bk} + A \sum \fg^{i\bi}
 \end{aligned}
\end{equation}
where $\lambda_i = (n-1) (\mu)_i - \mu T_{ji}^j$ and
\begin{equation}
\label{gblq-C151}
 A = |\mu| \sup_k |R_{k\bk} - S_{k\bk}| + |\mu| |T|^2 + |\nabla (\mu)| |T|
  + \Delta (\mu) - n \inf_k (\mu)_{k\bk}.
 \end{equation}
\end{lemma}

\begin{proof}
Since $e^\phi (n \mu + \Delta u)$ achieves its maximum at $p$,
\begin{equation}
\label{gblq-C80}
  \frac{(n \mu + \Delta u )_i}{n \mu + \Delta u} + \phi_i = 0, \;\;
\frac{(n \mu + \Delta u)_{\bi}}{n \mu + \Delta u} + \phi_{\bi} = 0,
\end{equation}
\begin{equation}
\label{gblq-C90}
 \frac{(n \mu + \Delta u)_{i\bi}}{n \mu+ \Delta u}
-  \frac{|(n \mu + \Delta u)_{i}|^2}{(n \mu + \Delta u)^2}
  + \phi_{i\bi} \leq 0.
\end{equation}
Note that
\[ (n \mu + \Delta u)_{i} = \sum_j (u_{i\bj} + \mu g_{i\bj})_j + \lambda_i\]
by \eqref{gblq-C30} and \eqref{coord}.
We have by \eqref{gblq-C80},
\begin{equation}
\label{gblq-C130}
\begin{aligned}
|(n \mu + \Delta u)_{i}|^2 \,
= \, & \sum_j |(u_{i\bj} + \mu g_{i\bj})_j|^2
+ 2 \fRe \{(n \mu + \Delta u)_{i} \bar{\lambda_{i}}\} - |\lambda_i|^2 \\
= \sum_j |(u_{i\bj} \, &  + \mu g_{i\bj})_j|^2
  - 2 (n \mu + \Delta u) \fRe\{\phi_{i} \bar{\lambda_{i}}\} - |\lambda_i|^2.
 \end{aligned}
\end{equation}
By Cauchy-Schwarz inequality,
\begin{equation}
\label{gblq-C140}
  \begin{aligned}
 \sum_{i,j}  \fg^{i\bi} |(u_{i\bj} + \mu g_{i\bj})_j|^2
   = & \; \sum_i \fg^{i\bi} \Big|\sum_j
     \fg_{j\bj}^{1/2} \; \fg_{j\bj}^{- 1/2} \;
               (u_{i\bj} + \mu g_{i\bj})_j\Big|^2 \\
\leq   & \; (n \mu + \Delta u)  \fg^{i\bi} \fg^{j\bj}
      |u_{i\bj j} + (\mu g_{i\bj})_j|^2.
 \end{aligned}
\end{equation}
From \eqref{gblq-C90}, \eqref{gblq-C120}, \eqref{gblq-C130} and
\eqref{gblq-C140} we derive \eqref{gblq-C150}.
\end{proof}

Let $\phi = e^{\eta (u)}$ with $\eta \geq 0$, $\mu \eta' < 0$, and $\eta'' \geq 0$. We have
\begin{equation}
\label{gblq-C155}
  \phi_{i} =   e^\eta \eta' u_i, \;\;
  \phi_{i\bi} = e^\eta [\eta' u_{i\bi} + (\eta'' + \eta'^2) |u_i|^2].
 \end{equation}
Therefore,
\begin{equation}
\label{gblq-C160}
  \begin{aligned}
  2 \fg^{i\bi} \fRe\{\phi_{i} \bar{\lambda_{i}}\}
    = \,& 2 e^\eta \eta' \fg^{i\bi} \fRe\{u_{i} \bar{\lambda_{i}}\}
 \geq - e^\eta \fg^{i\bi} (|\lambda_{i}|^2 + \eta'^2 |u_{i}|^2).
 \end{aligned}
\end{equation}

Suppose now that both $\psi$ and $\mu$ are independent of $u$.
Plugging \eqref{gblq-C155} and \eqref{gblq-C160} into \eqref{gblq-C150},
we see that
\begin{equation}
\label{gblq-C170}
  \begin{aligned}
0  \geq \,
  & (n \mu + \Delta u) \Big(e^{-\eta} \inf_{j,k} R_{j\bj k\bk}
    - \mu \eta'\Big) \sum \fg^{i\bi} \\
  & + \eta'^2 (n \mu + \Delta u - 1) \sum \fg^{i\bi} |u_{i}|^2
    + n \eta' (n \mu + \Delta u) \\
  & - n^2 e^{-\eta} \inf_{j,k} R_{j\bj k\bk} + e^{-\eta} \Delta f - C_3 \sum \fg^{i\bi}
\end{aligned}
\end{equation}
where 
\begin{equation}
\label{gblq-C180}
   \begin{aligned}
 C_3 = A + n^2 |\nabla \mu|^2 + n^2 \mu^2 |T|^2
\end{aligned}
\end{equation}
and $A$ is given in \eqref{gblq-C151}.

Following Yau~\cite{Yau78} we shall make use of the inequality
\begin{equation}
\label{gblq-C190}
 \Big(\sum \fg^{i\bi}\Big)^{n-1}
\geq \frac{\sum \fg_{i\bi}}{\fg_{1\bar{1}} \cdots \fg_{n\bar{n}}}
   = \frac{n \mu + \Delta u}{\det (\mu g_{i\bj} + u_{i\bj})}
   = \frac{n \mu + \Delta u}{\psi}.
  \end{equation}
Choosing $\eta = A (U - \nu u)$ where $\nu = \mu/|\mu|$,
$U = \sup_{M} \nu u$ and $A > 0$ is a constant such that
\begin{equation}
\label{gblq-C200}
 |\mu| A + e^{-\eta} \inf_{j,k} R_{j\bj k\bk} \geq 2,
\end{equation}
we see from \eqref{gblq-C170} and \eqref{gblq-C190} that
\begin{equation}
\label{gblq-C210}
 (n \mu + \Delta u)^{\frac{n}{n-1}}
  + n \nu A \psi^{\frac{1}{n-1}} (n \mu + \Delta u)
    + \psi^{\frac{1}{n-1}} \Delta f - n^2 e^{-\eta} \psi^{\frac{1}{n-1}} \inf_{j,k} R_{j\bj k\bk} \leq 0
\end{equation}
provided that
\begin{equation}
\label{gblq-C220}
 n \mu + \Delta u \geq 1 + C_3.
\end{equation}
This gives us a bound $(n \mu + \Delta u) (0) \leq C$ which depends
on $|u|_{C^0 (\bM)}$, $|\psi^{\frac{1}{n-1}}|_{C^2 (M)}$, $|\mu|_{C^2 (M)}$ and
geometric quantities of $(M, g)$. Finally,
\begin{equation}
\label{gblq-C230}
 \sup_M (n \mu + \Delta u)
      \leq e^{\phi (0)- \inf_M \phi} (n \mu + \Delta u) (0) \leq C.
\end{equation}

We have therefore proved the following.

\begin{proposition}
\label{gblq-prop-C10}
Suppose both $\mu$ and $\psi$ are independent of $u$.
Then
\begin{equation}
\label{gblq-C240}
\max_{\bM} \Delta u \leq
 C (1 + \max_{\partial M}  \Delta u)
\end{equation}
where $C > 0$ depends on
\[ |u|_{C^0 (\bM)}, 
  \; |\psi^{\frac{1}{n-1}}|_{C^2 (\bM)}, \; |\mu|_{C^2 (\bM)},
  \; \sup_{M} \frac{1}{|\mu|},
\; \inf_{j, k} R_{j\bj k\bk}, \; \sup_k |R_{k\bk} - S_{k\bk}|, \; |T|^2. \]
\end{proposition}

If $\psi$ and $\mu$ depend also on $u$ the estimate \eqref{gblq-C240} still
holds with $C$ depending in addition on $\sup_M |\nabla u|$. Indeed,
in places of \eqref{gblq-C170} we have
\begin{equation}
\label{gblq-C170'}
  \begin{aligned}
0  \geq \,
  & (n \mu + \Delta u) \Big(e^{-\eta} \inf_{j,k} R_{j\bj k\bk}
     - (n+1) e^{-\eta} |\mu_u| - \mu \eta'\Big) \sum \fg^{i\bi} \\
 &   + (n \eta' + e^{-\eta} f_u)  (n \mu + \Delta u) \\
  & - n^2 e^{-\eta} \inf_{j,k} R_{j\bj k\bk} - C'_3 \sum \fg^{i\bi}
\end{aligned}
\end{equation}
where $C'_3$ depends on $C_3$, $|u|_{C^1 (\bM)}$, as well as
the derivatives of $\psi^{\frac{1}{n-1}}$ and $\mu$. This again will give a bound for
$(n \mu + \Delta u) (0)$ and therefore \eqref{gblq-C230}.

\bigskip

\section{Boundary estimates for second derivatives}
\label{gblq-B}
\setcounter{equation}{0}
\medskip

In this section we derive
{\em a priori} estimates for second derivatives (the {\em real} Hessian) on the boundary
\begin{equation}
 \label{cma-37}
\max_{\partial M} |\nabla^2 u| \leq C.
\end{equation}
In order to track the dependence on the curvature and torsion of the
estimates we shall use covariant derivatives in this section. So we begin
with a brief review of formulas for changing the orders of covariant
derivatives which we shall also need in Section~\ref{gblq-R}.

In local coordinates $z = (z_1, \ldots, z_n)$,
$z_j = x_j + \sqrt{-1} y_j$, we shall use notations such as
\[ v_i = \nabla_{\frac{\partial}{\partial z_i}} v, \;
   v_{ij} = \nabla_{\frac{\partial}{\partial z_j}} \nabla_{\frac{\partial}{\partial z_i}} v, \;
   v_{x_i} = \nabla_{\frac{\partial}{\partial x_i}} v, \; \mbox{etc.} \]
Recall that
\begin{equation}
\label{gblq-B105}
 \begin{aligned}
& v_{i\bj} - v_{\bj i}  = 0, \;\;
  v_{ij} - v_{ji}  = T_{ij}^l v_l.
\end{aligned}
\end{equation}
By straightforward calculations,
\begin{equation}
\label{gblq-B145}
\left\{
\begin{aligned}
v_{i \bj \bk} - v_{i \bk \bj} = \, & \ol{T_{jk}^l} v_{i\bl},  \\
 v_{i \bj k} - v_{i k \bj} = \, & - g^{l\bm} R_{k \bj i \bm} v_l, \\
 v_{i j k} - v_{i k j} = \, & g^{l\bm} R_{jki\bm} v_l + T_{jk}^l v_{il}.
\end{aligned}
 \right.
\end{equation}
Therefore,
\begin{equation}
\label{gblq-B150}
 \begin{aligned}
 v_{i \bj k} - v_{k i \bj}
  = \,& (v_{i \bj k} - v_{i k \bj}) + (v_{i k \bj} - v_{k i \bj}) \\
  = \,& - g^{l\bm} R_{k \bj i \bm} v_l + T_{ik}^l v_{l\bj}
         + \nabla_{\bj} T_{ik}^l v_l \\
  = \,& - g^{l\bm} R_{i \bj k \bm} v_l + T_{ik}^l v_{l\bj}
\end{aligned}
\end{equation}
by \eqref{cma-K115},
and
\begin{equation}
\label{gblq-B155}
 \begin{aligned}
 v_{ijk} - v_{kij}
   =  \,& (v_{ijk} - v_{ikj}) + (v_{ikj} - v_{kij}) \\
   =  \,&  g^{l\bm} R_{jk i \bm} v_l + T_{jk}^l v_{il}
          + T_{ik}^l  v_{l j}  + \nabla_{j} T_{ik}^l v_l.
\end{aligned}
\end{equation}

Since
\[ \frac{\partial}{\partial x_k}
     = \frac{\partial}{\partial z_k} + \frac{\partial}{\partial \bz_k}, \;\;
   \frac{\partial}{\partial y_k}
 = \sqrt{-1} \Big(\frac{\partial}{\partial z_k}
  - \frac{\partial}{\partial \bz_k}\Big), \]
we see that
\begin{equation}
\label{gblq-B110}
\left\{
\begin{aligned}
v_{z_i x_j} -  v_{x_j z_i} = \,&  v_{ij} - v_{ji}  = T_{ij}^l v_l,  \\
v_{z_i y_j} -  v_{y_j z_i} = \,& \sqrt{-1} (v_{ij} - v_{ji})
                           = \sqrt{-1} T_{ij}^l v_l,
\end{aligned} \right.
\end{equation}
\begin{equation}
\label{gblq-B120}
\begin{aligned}
v_{z_i \bz_j x_k} - v_{x_k z_i \bz_j}
     = \,&  (v_{i\bj k} + v_{i\bj \bk}) - (v_{ki\bj} + v_{\bk i \bj}) \\
     = \,&  (v_{i\bj k} - v_{ki\bj}) + (v_{i\bj \bk} - v_{i \bk \bj}) \\
    = \,& - g^{l\bm} R_{i \bj k \bm} v_l + T_{ik}^l v_{l\bj}
          + \ol{T_{jk}^l} v_{i\bl},
  \end{aligned} 
\end{equation}
and, similarly,
\begin{equation}
\label{gblq-B130}
\begin{aligned}
v_{z_i z_{\bj} y_k} -  v_{y_k z_i z_{\bj}}
    = \,& \sqrt{-1} ((v_{i\bj k} - v_{i\bj \bk})
                      - (v_{ki\bj} -  v_{\bk i \bj})) \\
    = \,& \sqrt{-1} ((v_{i\bj k} - v_{ki\bj})
                      - (v_{i\bj \bk} - v_{i \bk \bj})) \\
    = \,&  \sqrt{-1} (- g^{l\bm} R_{i \bj k \bm} v_l + T_{ik}^l v_{l\bj}
                       - \ol{T_{jk}^l} v_{i\bl}).
\end{aligned} 
\end{equation}

For convenience we set
\[ 
t_{2k-1} = x_k,\;  t_{2k} = y_k, \, 1 \leq k \leq n-1; \;
  t_{2n-1} = y_n, \; t_{2n} = x_n. \]
By \eqref{gblq-B120}, \eqref{gblq-B130} and the identity
 \begin{equation}
\label{gblq-B170'}
 \fg^{i\bj} T^l_{ki} u_{l\bj} = T^i_{ki} - \mu \fg^{i\bj} T^l_{ki} g_{l\bj}
 \end{equation}
we obtain for all $1 \leq \alpha \leq 2n$,
\begin{equation}
\label{gblq-B50}
\begin{aligned}
|\fg^{i\bj}  (u_{t_{\alpha} i\bj} - u_{i\bj t_{\alpha}})|
 = \, & |\fg^{i\bj} u_{t_{\alpha} i\bj} - (f)_{t_{\alpha}}| \\
 \leq \, & 2 |T| + (|\mu| |T| + |R| + |\nabla u| |\nabla T|)
     \sum \fg^{i\bj} g_{i\bj}.
\end{aligned}
\end{equation}
We also record here the following identity which we shall use later:
for a function $\eta$,
\begin{equation}
\label{gblq-B190}
 \begin{aligned}
  \fg^{i\bj} \eta_i u_{x_n \bj}
   = &\, \fg^{i\bj} \eta_i (2 u_{n \bj} + \sqrt{-1} u_{y_n \bj}) \\
   = &\, 2 \eta_n - 2 \mu \fg^{i\bj} \eta_i g_{n \bj}
          + \sqrt{-1} \fg^{i\bj}  \eta_i u_{y_n \bj}.
 \end{aligned}
 \end{equation}

We now start to derive \eqref{cma-37}.
We assume
\begin{equation}
\label{cma-38}
|u| + | \nabla u| \leq K \;\; \mbox{in $\bar{M}$}.
\end{equation}
Set
\begin{equation}
\label{cma-39}
 \ul{\psi} \equiv
       \min_{|u| \leq K, \, z \in \bar{M}} \psi (z, u) > 0, \;\;
   \bar{\psi} \equiv
       \max_{|u| \leq K, \, z \in \bar{M}} \psi (z, u).
\end{equation}

\vspace{.2in}
Let $\sigma$ be the distance function to $\partial M$. Note that
$|\nabla \sigma| = \frac{1}{2}$ on $\partial M$.
 There exists
$\delta_0 > 0$ such that $\sigma$ is smooth and $\nabla \sigma \neq 0$
in
\[ M_{\delta_0} := \{z \in M: \sigma (z) < \delta_0\}, \]
which we call the $\delta_0$-neighborhood of $\partial M$.
We can therefore write
\begin{equation}
\label{gblq-B10}
 u - \ul{u} = h \sigma, \;\; \mbox{in $M_{\delta_0}$}
\end{equation}
where $h$ is a smooth function.

Consider a boundary point $p \in \partial M$. We choose local
coordinates $z = (z_1, \ldots, z_n)$, $z_j = x_j + i y_j$,
around $p$ in a neighborhood which we assume to be contained in
$M_{\delta_0}$ such
$\frac{\partial}{\partial x_n}$ is the interior normal direction to
$\partial M$ at $p$ where we also assume $g_{i\bj} = \delta_{ij}$;
(Here and in what follows we identify $p$ with $z = 0$.) for later
reference we call such local coordinates {\em regular} coordinate
charts.

By \eqref{gblq-B10} we have
\[  (u - \ul{u})_{x_n} = h_{x_n} \sigma + h \sigma_{x_n} \]
and
\[ (u - \ul{u})_{j\bk} =h_{j\bk} \sigma + h \sigma_{j\bk}
                      + 2 \, \fRe \{h_j \sigma_{\bk}\}. \]
Since $\sigma = 0$ on $\partial M$ and
$\sigma_{x_n} (0) = 2 |\nabla \sigma| = 1$, we see that
\[ (u - \ul{u})_{x_n} (0) = h (0) \]
and
\begin{equation}
\label{cma-60'}
(u-\ul{u})_{j\bk}(0) = (u-\ul{u})_{x_n} (0) \sigma_{j\bk}(0)
\;\; j, k < n.
\end{equation}
Similarly,
\begin{equation}
\label{cma-60}
(u-\ul{u})_{t_{\alpha}t_{\beta}} (0) = - (u-\ul{u})_{x_n} (0)
 \sigma_{t_{\alpha} t_{\beta}},
\;\; \alpha, \beta < 2n.
\end{equation}
It follows that
\begin{equation}
|u_{t_{\alpha} t_{\beta}}(0)| \leq C, \;\;\;\; \alpha, \beta < 2n
\label{cma-70}
\end{equation}
where $C$ depends on $|u|_{C^1 (\bM)}$, $|\ul{u}|_{C^1 (\bM)}$, and
the principal curvatures of $\partial M$.

\vspace{.2in}
To estimate $u_{t_{\alpha} x_n} (0)$ for $\alpha \leq 2n$,
we will follow \cite{Guan98b} and employ a barrier function of the form
\begin{equation}
\label{cma-E85}
v = (u - \ul{u}) + t \sigma - N \sigma^2,
\end{equation}
where $t, N$ are positive constants to be determined.
Recall that $\ul{u} \in C^2$ and $\omega_{\ul{u}} > 0$ in a neighborhood of
$\partial M$. We may assume that there exists $\epsilon > 0$ such that
$\omega_{\ul{u}} > \epsilon \omega$ in $M_{\delta_0}$. Locally, this gives
\begin{equation}
\label{cma-E86}
\{\ul{u}_{j\bk} + \mu g_{j\bk}\} \geq \epsilon \{g_{j\bk}\}.
\end{equation}
The following is the key ingredient in our argument.

\begin{lemma}
\label{cma-lemma-20}
For $N$ sufficiently large and $t, \delta$ sufficiently small,
\[ \begin{aligned}
  \fg^{i\bj} v_{i\bj}
    \leq - \frac{\epsilon}{4} &\,\Big(1 + \sum \fg^{i\bj} g_{i\bj}\Big)
   \;\;\; \mbox{in} \;\; \Omega_{\delta}, \\
  v & \, \geq 0 \;\; \mbox{on} \;\; \partial \Omega_{\delta}
  \end{aligned} \]
where $\Omega_{\delta} = M \cap B_{\delta}$ and $B_{\delta}$ is the
(geodesic) ball of radius $\delta$ centered at $p$.
\end{lemma}

\begin{proof}
This lemma was first proved in \cite{Guan98b} for domains in $\bfC^n$.
For completeness we include the proof here with minor modifications.
By \eqref{cma-E86} we have
\begin{equation}
\label{cma-E90}
\fg^{i\bj} (u_{i\bj} -\ul{u}_{i\bj})
   = \fg^{i\bj} (u_{i\bj} + \mu g_{i\bj} -\ul{u}_{i\bj} - \mu g_{i\bj})
 \leq n - \epsilon \sum \fg^{i\bj} g_{i\bj}.
\end{equation}
Obviously,
\[ \fg^{i\bj} \sigma_{i\bj} \leq C_1 \sum \fg^{i\bj}  g_{i\bj} \]
for some constant $C_1 > 0$ under control.
Thus
\[  \fg^{i\bj} v_{i\bj} \leq n
    + \{C_1 (t + N \sigma) - \epsilon\} \sum \fg^{i\bj} g_{i\bj}
- 2 N \fg^{i\bj} \sigma_i \sigma_{\bj} \;\;\; \mbox{in} \;\;
\Omega_{\delta}.\]

Let $\lambda_1 \leq \cdots \leq \lambda_n$ be the eigenvalues
of $\{u_{i\bj} + \mu g_{i\bj}\}$ (with respect to $\{g_{i\bj}\}$).
We have $\sum \fg^{i\bj}  g_{i\bj} = \sum \lambda_k^{-1}$
and
\begin{equation}
\label{cma-E95}
\fg^{i\bj} \sigma_i \sigma_{\bj} \geq  \frac{1}{2 \lambda_n}
\end{equation}
since $|\nabla \sigma| \equiv \frac{1}{2}$ where $\sigma$ is smooth.
By the arithmetic-geometric mean-value inequality,
\[ \frac{\epsilon}{4}
   \sum \fg^{i\bj}  g_{i\bj} + \frac{N}{\lambda_n}
\geq
  \frac{n \epsilon}{4} (N \lambda_1^{-1} \cdots \lambda_n^{-1})^{\frac{1}{n}}
  \geq \frac{n \epsilon  N^{\frac{1}{n}}}{4 \psi^{\frac{1}{n}}}
  \geq c_1 N^{\frac{1}{n}} \]
for some constant $c_1 > 0$ depending on the upper bound of $\psi$.

We now fix $t > 0$ sufficiently small and $N$ large so that
 $c_1 N^{1/n} \geq 1 + n + \epsilon$ and $C_1 t \leq \frac{\epsilon}{4}$.
Consequently,
\[  \fg^{i\bj} v_{i\bj}
    \leq - \frac{\epsilon}{4} \Big(1 + \sum \fg^{i\bj} g_{i\bj}\Big)
   \;\;\; \mbox{in} \;\; \Omega_{\delta} \]
if we require $\delta$ to satisfy
$C_1 N \delta \leq \frac{\epsilon}{4}$ in $\Omega_{\delta}$.

On $\partial M \cap B_{\delta}$ we have $v = 0$.
On $M \cap \partial B_{\delta}$,
\[ v \geq t \sigma - N \sigma^2
     \geq  (t - N \delta) \sigma \geq 0 \]
if we require, in addition, $N \delta \leq t$.
\end{proof}

\begin{remark}
For the real Monge-Amp\`ere equations,
Lemma~\ref{cma-lemma-20} was proved in \cite{Guan98a} both for domains
in $\bfR^n$ and in general Riemannian manifolds, improving earlier results
in \cite{HRS}, \cite{GS93} and \cite{GL96}.
\end{remark}

\begin{lemma}
\label{cma-lemma-30}
Let $w \in C^2 (\ol{\Omega_{\delta}})$.
Suppose that $w$ satisfies
\[ \fg^{i\bj} w_{i\bj} \geq - C_1 \Bigl(1 + \sum \fg^{i\bj} g_{i\bj}\Big) \;\;
\mbox{in $\Omega_{\delta}$} \]
and
\[ w \leq C_0 \rho^2 \;\; \mbox{on $ B_{\delta} \cap \partial M$},
   \;\; w (0) = 0\]
where $\rho$ is the distance function to the point $p$ (where $z = 0$) on
$\partial M$. Then
$ w_{\nu} (0) \leq C$,
where $\nu$ is the interior unit normal to $\partial M$, and
$C$ depends on $\epsilon^{-1}$, $C_0$, $C_1$,
$|w|_{C^0 (\ol{\Omega_{\delta}})}$,
$|u|_{C^1 (\bar{M})}$ and the constants $N$, $t$ and $\delta$ determined
in Lemma~\ref{cma-lemma-20}.
\end{lemma}

\begin{proof}
By Lemma~\ref{cma-lemma-20},
$A v + B \rho^2 - w \geq 0$ on $\partial \Omega_{\delta}$
and
\[ \fg^{i\bj} (A v + B \rho^2 - w)_{i\bj} \leq 0 \;\;
     \mbox{in $\Omega_{\delta}$} \]
when $A \gg B$ and both are sufficiently large. By the maximum principle,
\[ A v + B \rho^2 - w \geq 0 \;\; \mbox{in $\ol{\Omega_{\delta}}$}. \]
Consequently,
\[ A v_{\nu} (0) - w_{\nu} (0) = D_{\nu} (A v + B \rho^2 - w) (0) \geq 0 \]
since $A v + B \rho^2 - w = 0$ at the origin.
\end{proof}

\vspace{.2in}
We next apply Lemma~\ref{cma-lemma-30} to estimate
$u_{t_{\alpha} x_n} (0)$ for $\alpha < 2n$. For fixed $\alpha < 2n$,
we write $\eta = \sigma_{t_{\alpha}}/\sigma_{x_n}$ and
define
\[ \mathcal{T} = \nabla_{\frac{\partial}{\partial t_{\alpha}}}
     - \eta \nabla_{\frac{\partial}{\partial x_n}}. \]
We wish to apply Lemma~\ref{cma-lemma-30}  to
\[ w = (u_{y_n} - \varphi_{y_n})^2 \pm  \cT (u - \varphi). \]

By (\ref{cma-38}),
\[ | \mathcal{T} (u - \varphi)| + (u_{y_n} - \varphi_{y_n})^2 \leq C
    \;\; \mbox{in} \; \Omega_{\delta}. \]
On $\partial M$ since $u - \varphi = 0$ and
$\cT$ is a tangential differential operator, we have
\[  \mathcal{T} (u - \varphi) = 0 \;\; \mbox{on}\;\partial M \cap B_{\delta} \]
and, similarly,
\begin{equation}
\label{cma-175}
(u_{y_n} - \varphi_{y_n})^2 \leq C \rho^2
  \;\; \mbox{on} \; \partial M \cap B_{\delta}.
\end{equation}

We compute next
\begin{equation}
\label{gblq-B160}
\fg^{i\bj}  (\mathcal{T} u)_{i\bj}
   = \fg^{i\bj}  (u_{t_\alpha i\bj} + \eta u_{x_n i\bj})
 + \fg^{i\bj} \eta_{i\bj} u_{x_n} + 2 \fg^{i\bj}  \fRe \{\eta_i u_{x_n \bj}\}.
\end{equation}

By \eqref{gblq-B50} and \eqref{gblq-B190},
\begin{equation}
\label{gblq-B180}
 |\fg^{i\bj}  (u_{t_\alpha i\bj} + \eta u_{x_n i\bj})|
   \leq | \mathcal{T} (f) |
   + C_1 (|T| + |R| + |\nabla T|) \Big(1 + \sum \fg^{i\bj}  g_{i\bj}\Big)
\end{equation}
and
\begin{equation}
\label{gblq-B200}
2 |\fg^{i\bj} \fRe\{\eta_i u_{x_n \bj}\}| \leq \fg^{i\bj} u_{y_n i} u_{y_n \bj}
+ C_2 \Big(1 + \sum \fg^{i\bj}  g_{i\bj}\Big)
\end{equation}
where $C_1$ and $C_2$ are independent of the curvature and torsion.
 Applying \eqref{gblq-B50} again, we derive
\begin{equation}
\label{gblq-B210}
\begin{aligned}
\fg^{i\bj} [(u_{y_n} - \varphi_{y_n})^2]_{i\bj}
   = \,& 2 \fg^{i\bj}
         (u_{y_n} - \varphi_{y_n})_i (u_{y_n} - \varphi_{y_n})_{\bj} \\
       & + 2 (u_{y_n} - \varphi_{y_n})
         \fg^{i\bj} (u_{y_n} -  \varphi_{y_n})_{i\bj} \\
\geq \,& \fg^{i\bj} u_{y_n i} u_{y_n \bj}
          - 2 \fg^{i\bj} \varphi_{y_n i} \varphi_{y_n \bj} \\
       & + 2 (u_{y_n} - \varphi_{y_n})
          \fg^{i\bj} (u_{y_n i\bj} - \varphi_{y_n i\bj}) \\
\geq \,& \fg^{i\bj} u_{y_n i} u_{y_n \bj}
          - |(f)_{y_n}| - C_3 \\
       & - C_4 (1 + |T| + |R| + |\nabla T|) \sum \fg^{i\bj}  g_{i\bj}.
 \end{aligned}
 \end{equation}
Finally, combining \eqref{gblq-B160}-\eqref{gblq-B210} we obtain
\begin{equation}
\label{gblq-B230}
 \fg^{i\bj} [(u_{y_n} - \varphi_{y_n})^2 \pm \mathcal{T} (u - \varphi)]_{i\bj}
  \geq - C \Big(1 + |D f| + \sum \fg^{i\bj}  g_{i\bj}\Big)
  \;\; \mbox{in} \;\Omega_{\delta}
 \end{equation}
 where $C = C_0 (1 + |R| + |T| + |\nabla T|)$ with $C_0$ independent of
 the curvature and torsion.

Consequently, we may apply Lemma~\ref{cma-lemma-30} to
$w = (u_{y_n} - \varphi_{y_n})^2 \pm  \cT (u - \varphi)$ to obtain
\begin{equation}
\label{cma-180}
|u_{t_{\alpha} x_n} (0)| 
                         \leq C, \;\;\;\; \alpha < 2n.
\end{equation}
By \eqref{gblq-B110} we also have
\begin{equation}
\label{cma-180'}
|u_{x_n t_{\alpha}} (0)| 
                         \leq C, \;\; \alpha < 2n.
\end{equation}

\vspace{.2in}
It remains to establish the estimate
\begin{equation}
|u_{x_n x_n} (0)| \leq  C.
\end{equation}
Since we have already derived
\begin{equation}
\label{cma-190}
|u_{t_{\alpha} t_{\beta}} (0)|, \; |u_{t_{\alpha} x_n} (0)|, \;
|u_{x_n t_{\alpha}} (0)| \leq C,
\;\;\;\; \alpha, \beta < 2n,
\end{equation}
it suffices to prove
\begin{equation}
\label{cma-200}
0 \leq \mu + u_{n \bn} (0) =  \mu + u_{x_n x_n} (0) + u_{y_n y_n} (0) \leq C.
\end{equation}

Expanding $\det (u_{i\bj} + \mu g_{i\bj})$, we have
\begin{equation}
\label{gblq-B310}
 \det (u_{i\bj} (0) + \mu g_{i\bj}) = a (u_{n\bn} (0) + \mu) + b
\end{equation}
where
\[ a = \det (u_{\alpha \bar{\beta}} (0)+ \mu g_{\alpha \bar{\beta}})
   |_{\{1 \leq \alpha, \beta \leq n-1\}} \]
and $b$ is bounded in view of \eqref{cma-190}.
Since $\det (u_{i\bj} + \mu g_{i\bj})$ is bounded,
 we only have to derive an {\em a priori} positive lower bound for $a$, which
is equivalent to
\begin{equation}
\label{cma-80}
\sum_{\alpha, \beta < n}
  u_{{\alpha} \bar{\beta}} (0) \xi_{\alpha} \bar{\xi}_{\beta} \geq c_0 |\xi|^2, \;\;
\;\; \forall \, \xi 
\in {\bfC}^{n-1}
\end{equation}
for a uniform constant $ c_0 > 0$.

\begin{proposition}
\label{prop-cma-10}
 There exists $c_0 = c_0 (\ul{\psi}^{-1},
\varphi, \ul{u}) > 0$ such that (\ref{cma-80})
holds.
\end{proposition}

\begin{proof}
Let $T_C \partial M \subset T_C M $
be the complex tangent bundle of $\partial M$ and
\[ T^{1,0} \partial M =  T^{1,0} M \cap T_C \partial M
   = \Big\{\xi \in T^{1,0} M: d \sigma (\xi) = 0\Big\}. \]
 In local coordinates,
\[ T^{1,0} \partial M
   = \Big\{\xi = \xi_i \frac{\partial}{\partial z_i} \in T^{1,0} M:
        \sum \xi_i \sigma_i = 0 \Big\}. \]
It is enough to establish a positive lower bound for
\[ m_0 = 
       \min_{\xi \in T^{1,0} \partial M,  |\xi| = 1}
      \omega_u (\xi, \bar{\xi}). \]

We assume that $m_0$ is attained at a point $p \in \partial M$ and
choose regular local coordinates 
around $p$
as before such that
\[ m_0 = \omega_u \Big(\frac{\partial}{\partial z_1}, \frac{\partial}{\partial \bz_1}\Big)
       = u_{1\bar{1}} (0) + \mu. \]
One needs to show
\begin{equation}
\label{cma-90}
m_0 = u_{1\bar{1}} (0) + \mu \geq c_0 > 0.
\end{equation}
By \eqref{cma-60'},
\begin{equation}
\label{cma-100}
u_{1\bar{1}} (0) = \ul{u}_{1\bar{1}} (0)
        - (u-\ul{u})_{x_n} (0) \sigma_{1\bar{1}} (0).
\end{equation}
We can assume
$u_{1\bar{1}} (0) + \mu \leq \frac{1}{2} (\ul{u}_{1\bar{1}} (0) + \mu)$; otherwise we are done.
Thus 
\begin{equation}
\label{cma-100'}
(u-\ul{u})_{x_n} (0) \sigma_{1\bar{1}} (0) \geq \frac{1}{2} (\ul{u}_{1\bar{1}} (0) + \mu).
\end{equation}
It follows from \eqref{cma-38} that
\begin{equation}
\label{cma-101}
\sigma_{1\bar{1}} (0) \geq \frac{\ul{u}_{1\bar{1}} (0) + \mu}{2 C}
                      \geq \frac{\epsilon}{2 C} \equiv c_1 > 0
\end{equation}
where $C = \max_{\partial M} |\nabla (u-\ul{u})|$.

Let $\delta > 0$ be small enough so that
\[ \begin{aligned}
  w \equiv \, &  \Big|- \sigma_{z_n} \frac{\partial}{\partial z_1}
   + \sigma_{z_1} \frac{\partial}{\partial z_n}\Big|  \\
 = \, & \left(g_{1\bar{1}} |\sigma_{z_n}|^2 - 2 \fRe \{g_{1\bar{n}} \sigma_{z_n} \sigma_{\bz_1}\}
   + g_{n\bar{n}} |\sigma_{z_1}|^2\right)^{\frac{1}{2}} > 0 \;\;
\mbox{in $M \cap B_{\delta} (p)$}.
\end{aligned} \]
Define $\zeta = \sum \zeta_i \frac{\partial}{\partial z_i} \in T^{1,0} M$
in $M \cap B_{\delta} (p)$:
\[ \left\{ \begin{aligned}
   \zeta_1 & \, = - \frac{\sigma_{z_n}}{w}, \\
   \zeta_j & \, = 0, \; \; 2 \leq j \leq n-1, \\
   \zeta_n & \, = \frac{\sigma_{z_1}}{w}
   \end{aligned} \right.  \]
and
\[ \varPhi = (\varphi_{j\bk} + \mu g_{j\bk}) \zeta_j \bar{\zeta}_k
   - (u - \varphi)_{x_n} \sigma_{j\bk} \zeta_j \bar{\zeta}_k - u_{1\bar{1}}(0) - \mu. \]
Note that $\zeta \in T^{1,0} \partial M$ on $\partial M$ and $|\zeta| = 1$.
By \eqref{cma-60'},
\begin{equation}
\label{cma-103}
\varPhi = (u_{j\bk} + \mu g_{j\bk}) \zeta_j \bar{\zeta}_k - u_{1\bar{1}}(0) - \mu \geq 0
\;\;  \mbox{on $\partial M \cap B_{\delta} (p)$}
\end{equation}
and $ \varPhi (0) = 0$.

Write $G = \sigma _{i\bj} \zeta_i \bar{\zeta}_i$. We have
\begin{equation}
\label{gblq-B360}
 \begin{aligned}
\fg^{i\bj} \varPhi_{i\bj}
 \leq \, & - \fg^{i\bj} (u_{x_n} G)_{i\bj} + C \Big(1 + \sum u^{j\bj}\Big) \\
    = \, & - G \fg^{i\bj} u_{x_n i\bj} - 2 \fg^{i\bj} \fRe\{u_{x_n i} G_{\bj}\}
        + C \Big(1 + \sum \fg^{i\bj}  g_{i\bj}\Big) \\
 \leq \, & \fg^{i\bj} u_{y_n i} u_{y_n \bj}
       + C \Big(1 + \sum \fg^{i\bj}  g_{i\bj}\Big)
\end{aligned}
\end{equation}
by \eqref{gblq-B50} and \eqref{gblq-B190}.
It follows that
\begin{equation}
\label{cma-105}
\fg^{i\bj} [\varPhi - (u_{y_n} - \varphi_{y_n})^2]_{i\bj}
      \leq C \Big(1 + \sum \fg^{i\bj}  g_{i\bj}\Big)
        \;\; \mbox{in $M \cap B_{\delta} (p)$}.
\end{equation}
Moreover, by \eqref{cma-175} and \eqref{cma-103},
\[ (u_{y_n} - \varphi_{y_n})^2 - \varPhi  \leq C |z|^2 \;\;
 \mbox{on $\partial M \cap B_{\delta} (p)$}. \]
Consequently, we may apply Lemma~\ref{cma-lemma-30} to
\[ h = (u_{y_n} - \varphi_{y_n})^2 - \varPhi \]
to derive $\varPhi_{x_n} (0) \geq -C$
which, by \eqref{cma-101}, implies
\begin{equation}
\label{cma-310}
 u_{x_n x_n} (0) \leq \frac{C}{\sigma_{1\bar{1}} (0)} \leq \frac{C}{c_1}.
\end{equation}

In view of \eqref{cma-190} and \eqref{cma-310} we have an {\em a priori}
upper bound for all eigenvalues of $\{u_{i\bj} + \mu g_{i\bj}\}$ at $p$.
Since $\det (u_{i\bj} + \mu g_{i\bj}) \geq \ul{\psi} > 0$,
the eigenvalues of
$\{u_{i\bj} + \mu g_{i\bj}\}$ at $p$ must admit a positive lower bound,
i.e.,
\[ \min_{\xi \in T^{1,0}_p M, |\xi| = 1} (u_{i\bj} + \mu g_{i\bj}) \xi_i \bar{\xi}_j
    \geq c_0.\]
Therefore,
\[ m_0 = \min_{\xi \in T_p^{1,0} \partial M, |\xi| = 1}
         (u_{i\bj} + \mu g_{i\bj}) \xi_i \bar{\xi}_j
   \geq \min_{\xi \in T^{1,0}_p M, |\xi| = 1}
         (u_{i\bj} + \mu g_{i\bj}) \xi_i \bar{\xi_j}
    \geq c_0.\]
The proof of Proposition~\ref{prop-cma-10} is complete.
\end{proof}

We have therefore established (\ref{cma-37}).

\bigskip

\section{Estimates for the real Hessian and higher derivatives}
\label{gblq-R}
\setcounter{equation}{0}
\medskip

The primary goal of this section is to derive global estimates for the whole
(real) Hessian
\begin{equation}
\label{gblq-R5}
|\nabla^2 u| \leq  C  \;\; \mbox{on $\bM$}.
\end{equation}
This is equivalent to
\begin{equation}
\label{cma-410'}
 |u_{x_i x_j} (p)|, \; |u_{x_i y_j} (p)|, \;  |u_{y_i y_j} (p)| \leq C, \;\;
 \forall \, 1 \leq i, j \leq n
\end{equation}
in local coordinates  $z = (z_1, \ldots, z_n)$,
$z_j = x_j + \sqrt{-1} y_j$ with $g_{i\bj} (p) = \delta_{ij}$ for
any fixed point $p \in M$,
where the constant $C$ may depend on $|u|_{C^1 (M)}$,
$\sup_M \Delta u$, $\inf \psi > 0$, and the curvature and torsion of $M$ as
well as their derivatives. Once this is done we can apply the Evans-Krylov
Theorem to obtain global $C^{2, \alpha}$ estimates.

As in Section~\ref{gblq-B} we shall use covariant derivatives.
We start with communication formulas for the fourth order derivatives.
From  direct computation,
\begin{equation}
\label{gblq-R145}
\left\{
\begin{aligned}
v_{ij \bk \bl} - v_{ij \bl \bk} = \, & \ol{T_{kl}^q} v_{ij \bq}, \\
v_{ij k \bl} - v_{ij \bl k} = \, & g^{p\bq} R_{k\bl i\bq} v_{pj}
                                 + g^{p\bq} R_{k\bl j\bq} v_{ip}.
\end{aligned}
 \right.
\end{equation}
Therefore, by \eqref{gblq-B145}, \eqref{gblq-B150}, \eqref{gblq-B155},
\eqref{gblq-R145} and \eqref{cma-K115},
\begin{equation}
\label{gblq-R150}
 \begin{aligned}
v_{i \bj k \bl} - v_{k \bl i \bj}
 = \,& (v_{i \bj k \bl} - v_{k i \bj \bl})
         + (v_{k i \bj \bl} - v_{k i \bl \bj})
         + (v_{k i \bl \bj} - v_{k \bl i \bj}) \\
 = \,& \nabla_{\bl} (- g^{p\bq} R_{i \bj k \bq} v_p + T_{ik}^p v_{p\bj})
     + \ol{T_{jl}^q} v_{ki\bq} + g^{p\bq} \nabla_{\bj} (R_{i \bl k \bq} v_p)\\
 = \,& g^{p\bq} (R_{k\bl i\bq} v_{p\bj} - R_{i\bj k\bq} v_{p\bl})
       + T_{ik}^p v_{p\bj \bl} + \ol{T_{jl}^q} v_{ki\bq} \\
   & + g^{p\bq} (\nabla_{\bj} R_{i\bl k\bq} - \nabla_{\bl} R_{i\bj k\bq}) v_p
\end{aligned}
\end{equation}
and
\begin{equation}
\label{gblq-R155}
 \begin{aligned}
v_{i \bj k l} - v_{k l i \bj}
   = \,& v_{i \bj k l} - v_{k i \bj l}
         + v_{k i \bj l} - v_{k i l \bj}
         + v_{k i l \bj} - v_{k l i \bj} \\
   = \,& \nabla_{l} (- g^{p\bq} R_{i \bj k \bq} v_p + T_{ik}^p v_{p\bj})
     -  g^{p\bq} R_{l\bj k\bq} v_{pi} - g^{p\bq} R_{l\bj i\bq} v_{kp} \\
     & + \nabla_{\bj}  (g^{p\bq} R_{ilk\bq} v_p + T_{il}^p v_{kp})\\
   = \,& -g^{p \bq}R_{i\bj k \bq}v_{pl}-g^{p \bq}R_{i \bj l \bq}v_{kp}-g^{p \bq}R_{l \bj k \bq}v_{pi} \\
   & -g^{p \bq}[(\nabla_l R_{i \bj k \bq})+(\nabla_{\bj} R_{ilk \bq})]v_p\\
   & + [(\nabla_l T_{ik}^p)+g^{p \bq} R_{i l k \bq}]v_{p \bj}\\
   & + T_{ik}^p v_{p \bj l}+T_{il}^p v_{k p \bj}.
\end{aligned}
\end{equation}

Turning to the proof of \eqref{cma-410'}, it suffices to prove
the following.

\begin{proposition}
\label{gblq-prop-R10}
 There exists constant $C > 0$ depending on $|u|_{C^1 (\bM)}$,
$\sup_M \Delta u$ and $\inf \psi > 0$ such that
\begin{equation}
\label{gblq-R5'}
 \sup_{\tau \in T M, |\tau| = 1} u_{\tau \tau}  \leq C.
\end{equation}
\end{proposition}

\begin{proof}
Let
\[ N :=\sup_M \Big\{|\nabla u|^2 + A |\omega_u|^2
       + \sup_{\tau \in T M, |\tau| = 1} u_{\tau \tau} \Big\} \]
where $A$ is positive constant to be determined,
and assume that it is achieved at an interior point $p \in M$ and for some
unit vector $\tau \in T_p M$.
We choose local coordinates $z = (z_1, \ldots, z_n)$
such that $g_{i\bj} = \delta_{ij}$ and
$\{u_{i\bj}\}$ is diagonal at $p$. Thus $\tau$ can be written in the form
\[ \tau = a_j \frac{\partial}{\partial z_j} + b_j \frac{\partial}{\partial \bz_j},
\;\; a_j, b_j \in \bfC, \;\;  \sum a_j b_j = \frac{1}{2}.  \]

Let $\xi$ be a smooth unit vector field defined in a neighborhood of $p$
such that $\xi (p) = \tau$. Then the function
\[ Q = u_{\xi \xi} + |\nabla u|^2 + A |\omega_u|^2 \]
(defined in a neighborhood of $p$) attains it maximum at $p$ where,
therefore
\begin{equation}
\label{gblq-R20}
Q_i = u_{\tau \tau i} + u_k u_{i\bk} + u_{\bk} u_{ki}
   + 2 A (u_{k\bk}+\mu) (u_{k\bk i}+(\mu)_i) = 0 
\end{equation}
and
\begin{equation}
\label{gblq-R30}
\begin{aligned}
0 \geq \fg^{i\bi} Q_{i\bi}
 =  \,& \fg^{i\bi}(u_{k\bi} u_{i\bk} + u_{ki} u_{\bk \bi})
        +  \fg^{i\bi} (u_k u_{i\bk \bi} + u_{\bk} u_{ki\bi}) \\
      + \fg^{i\bi} \,& u_{\tau \tau i\bi}
        + 2 A \fg^{i\bi} (u_{k\bl i}+(\mu)_i)
    (u_{l\bk \bi}+(\mu)_{\bi}) \\
   &  + 2 A (u_{k\bk}+\mu) \fg^{i\bi} (u_{k\bk i\bi}+(\mu)_{i\bi}).
\end{aligned}
\end{equation}

Differentiating equation~(\ref{cma2-M10}) twice (using covariant derivatives),
by \eqref{gblq-B150} and \eqref{gblq-R150} we obtain
\begin{equation}
\label{gblq-R40}
 \fg^{i\bi} u_{ki\bi}
     = (f)_k + \fg^{i\bi} R_{i\bi k\bl} u_l
        - \fg^{i\bi}T_{ik}^l u_{l\bi} - (\mu)_k \sum \fg^{i\bi}
                         \geq (f)_k - C \sum \fg^{i\bi}
\end{equation}
and
\begin{equation}
\label{gblq-R50}
\begin{aligned}
\fg^{i\bi} u_{k\bk i\bi}
 \geq \,&  \fg^{i\bi} \fg^{j\bj} u_{i\bj k} u_{j\bi \bk}
    + \fg^{i\bi}  (T_{ik}^p u_{p\bi \bk} + \ol{T_{ik}^p} u_{ki\bp}) \\
     & + (f)_{k\bk} - C \sum \fg^{i\bi}
 \geq (f)_{k\bk} - C \Big(1 + \sum \fg^{i\bi}\Big).
\end{aligned}
\end{equation}
(Here we used \eqref{gblq-B150} again for the last inequality.)

Note that
\[ u_{\tau \tau i\bi} = a_k a_l u_{kl i\bi}
   + 2 a_k b_l  u_{k\bl i\bi} + b_k b_l  u_{\bk\bl i\bi}. \]
Using the formulas in
\eqref{gblq-R145}, \eqref{gblq-R150} and \eqref{gblq-R155} we obtain
\begin{equation}
\label{gblq-R60}
 \begin{aligned}
 \fg^{i\bi} u_{\tau \tau i\bi}
 \geq \,&  \fg^{i\bi} u_{i\bi \tau \tau} - C \fg^{i\bi} |T_{ik}^l u_{l\bi k}|
-  C \Big(1 + \sum_{k, l} |u_{kl}|\Big) \sum \fg^{i\bi} \\
 \geq \,& (f)_{\tau \tau} - C \fg^{i\bi} u_{l\bi k} u_{i\bl k}
-  C \Big(1 + \sum_{k, l} |u_{kl}|\Big) \sum \fg^{i\bi}.
\end{aligned}
\end{equation}
Plugging \eqref{gblq-R40}, \eqref{gblq-R50}, \eqref{gblq-R60} into
\eqref{gblq-R30} and using the inequality
\begin{equation}
\label{gblq-R70}
\begin{aligned}
 2 \fg^{i\bi} \,& (u_{k\bl i} + (\mu)_i)(u_{l\bk \bi}+(\mu)_{\bi}) \\
\,& \geq \fg^{i\bi} u_{k\bl i} u_{l\bk \bi} - \fg^{i\bi} (\mu)_i (\mu)_{\bi}
  \geq \fg^{i\bi} u_{k\bl i} u_{l\bk \bi} - C \sum \fg^{i\bi},
\end{aligned}
\end{equation}
we see that
\begin{equation}
\label{gblq-R80}
\fg^{i\bi} u_{ki} u_{\bk \bi} + (A-C) \fg^{i\bi} u_{k\bl i} u_{l\bk \bi}
- C \Big(1 + A +  \sum_{k, l} |u_{kl}|\Big) \Big(1 + \sum \fg^{i\bi}\Big)
 \leq 0.
\end{equation}
We now need the nondegeneracy of equation~\eqref{cma2-M10} which implies that
there is $\Lambda > 0$ depending on $\sup_M \Delta u$ and $\inf \psi > 0$
such that
\[ \Lambda^{-1} \{g_{i\bj}\} \leq \{\fg_{i\bj}\} \leq \Lambda \{g_{i\bj}\} \]
and, therefore,
\begin{equation}
\label{gblq-R90}
\left\{\begin{aligned}
& \sum \fg^{i\bi} \leq n \Lambda, \\
& \fg^{i\bi} u_{ki} u_{\bk \bi}
    \geq  \frac{1}{\Lambda} \sum_{i, k} | u_{ki}|^2.
\end{aligned} \right.
\end{equation}
Plugging these into \eqref{gblq-R80} and choosing $A$ large we derive
\[ \sum_{i, k} | u_{ki}|^2 \leq C.  \]
Consequently we have a bound $u_{\tau \tau} (p) \leq C$.
Finally, 
\[  \sup_{q \in M} \sup_{\tau \in T_q M, |\tau| = 1} u_{\tau \tau}
   \leq u_{\tau \tau} (p) + 2 \sup_M (|\nabla u|^2 + A |\omega_u|^2). \]
This completes the proof of \eqref{gblq-R5'}.
\end{proof}

We can now appeal to the Evans-Krylov Theorem (\cite{Evans}, \cite{Krylov82},
\cite{Krylov83}) for $C^{2, \alpha}$ estimates
\begin{equation}
\label{gblq-R100}
|u|_{C^{2, \alpha} (M)} \leq C.
\end{equation}
Higher order regularity and estimates now follow from the classical
Schauder theory for elliptic linear equations.

\begin{remark}
\label{gblq-remark-R10}
When $M$ is a K\"ahler manifold, Proposition~\ref{gblq-prop-R10} was recently
proved by Blocki~\cite{Blocki}. He observed that the estimate \eqref{gblq-R5'}
does not depend on $\inf \psi$ when $M$ has nonnegative bisectional curvature.
This is clearly also true in the Hermitian case.
\end{remark}

\begin{remark}
\label{gblq-remark-R20}
An alternative approach to the $C^{2, \alpha}$ estimate \eqref{gblq-R100} is
to use \eqref{gblq-I60} and the boundary estimate \eqref{cma-37} (in place of
\eqref{gblq-R5}) and apply an extension of the Evans-Krylov Theorem; see
Theorem 7.3, page 126 in \cite{CW} which only requires $C^{1, \alpha}$ bounds
for the solution. This was pointed out to us by Pengfei Guan to whom we wish
to express our gratitude.
\end{remark}

\bigskip

\section{$C^0$ estimates and existence}
\label{gblq-E}
\setcounter{equation}{0}
\medskip

In this section we complete the proof of
Theorem~\ref{gblq-th20}-\ref{gblq-th40} using the estimates established in
previous sections. We shall consider separately the Dirichlet problem and
the case of manifolds without boundary. In each case we need first to derive
$C^0$ estimates; the existence of solutions then can be proved by the continuity
method, possibly combined with degree arguments.

\subsection{Compact manifolds without boundary}
For the $C^0$ estimate on compact manifolds without boundary, we follow the
argument in \cite{Siu87}, \cite{Tian00} which simplifies the original proof
of Yau~\cite{Yau78}.

Let $M$ be a compact Hermitian manifold without boundary and $u$
 an admissible solution of equation~\eqref{cma2-M10},  $\sup_M u=-1$.
In this case we assume $\mu > 0$. (When $\mu < 0$ equation~\eqref{cma2-M10}
does not have solutions by the maximum principle.)
We write
\[ \chi = \sum_{k=0}^{n-1} (\m \o)^k \wedge (\o_u)^{n-1-k}. \]
Multiply the identity
$(\o_u)^n-(\m \o)^n = \frac{\sqrt{-1}}{2} \p \bar \p u \wedge \chi$
 by $(-u)^p$ and integrate over $M$,
\begin{equation}
\label{gblq-E120}
\begin{aligned}
 \int_M (-u)^p \,&[(\o_u)^n-(\m \o)^n]
   = \frac{\sqrt{-1}}{2} \int_M (-u)^p \p \bar \p u \wedge \chi \\
   = \,& \frac{p \sqrt{-1}}{2}
          \int_M (-u)^{p-1} \p u \wedge \bar \p u \wedge \chi
       + \frac{\sqrt{-1}}{2} \int_M (-u)^p \bar \p u \wedge \p \chi \\
   = \,& \frac{2p \sqrt{-1}}{(p+1)^2} \int_M
          \p (-u)^{\frac{p+1}{2}} \bar\p (-u)^{\frac{p+1}{2}} \wedge \chi
         - \frac{\sqrt{-1}}{2(p+1)} \int_M (-u)^{p+1} \p \bar \p \chi.
 \end{aligned}
\end{equation}
We now assume that $\p \bar\p (\mu \omega)^k = 0$, for all $1 \leq k \leq n-1$,
 which implies $\p \bar\p \chi = 0$, and that $\psi$ does not depends on $u$.
Since $\mu \omega > 0$ and $\omega_u \geq 0$,
we see that $(\mu \omega)^k \wedge (\omega_u)^{n-1-k} \geq 0$ for each $k$.
Therefore,
\begin{equation}
\label{gblq-E130}
\begin{aligned}
      \int_M  |\nabla (-u)^{\frac{p+1}{2}}|^2 \o^n
   = \,& 
       \frac{\sqrt{-1}}{2}
       \int_M \p (-u)^{\frac{p+1}{2}} \bar\p (-u)^{\frac{p+1}{2}}
               \wedge \omega^{n-1} \\
\leq \,& \frac{\sqrt{-1}}{2 \inf \mu^{n-1}}
         \int_M \p (-u)^{\frac{p+1}{2}} \bar\p (-u)^{\frac{p+1}{2}}
               \wedge \chi \\
   = \,& \frac{(p+1)^2}{2 p \inf \mu^{n-1}} \int_M (-u)^p (\psi-\mu^n)\o^n \\
 \leq \,& C \int_M (-u)^{p+1} \o^n.
 \end{aligned}
\end{equation}
After this we can derive a bound for $\inf u$ by
the Moser iteration method, following the argument in \cite{Tian00}.

If $\psi$ depends on $u$ and satisfies \eqref{gblq-I90}, a bound for
$\sup_M |u|$ follows directly from equation~\eqref{cma2-M10} by the maximum
principle.
Indeed, suppose $u (p) = \max_M u$ for some $p \in M$.
Then $\{u_{i\bj} (p)\} \leq 0$ and, therefore
\[ \mu^n \det g_{i\bj} \geq  \det (u_{i\bj} + \mu g_{ij})
= \psi (p, u(p)) \det g_{i\bj}. \]
This implies an upper bound $u (p) \leq C$ by \eqref{gblq-I90}.
That $\min_M u \geq -C$ follows from a similar argument.

\begin{proof}[Proof of Theorem~\ref{gblq-th30}]
We first consider the case that $\psi$ does not depend on $u$.
By assumption~\eqref{gblq-I70} we see that
\[ \int_M (\psi -1) \omega^n = 0 \]
is a necessary condition for the existence of admissible solutions,
and that the linearized operator,
$v \mapsto \fg^{i\bj} v_{i\bj}$, 
of equation~\eqref{gblq-I10} is self-adjoint.
So the continuity method proof in \cite{Yau78} works to give a unique
admissible solution $u \in \cH \cap C^{2, \alpha} (M)$ of \eqref{gblq-I10}
satisfying
\[ \int_M u \omega^n = 0. \]
The smoothness of $u$ follows from the Schauder regularity theory.

For the general case under the assumption $\psi_u \geq 0$, one
can still follow the proof of Yau~\cite{Yau78}. So we omit it here.
\end{proof}

\begin{proof}[Proof of Theorem~\ref{gblq-th40}]
The uniqueness follows easily from the assumption $\psi_u > 0$ and the
maximum principle. For the existence we make use of the continuity method.
For $0 \leq s \leq 1$ consider
\begin{equation}
\label{gblq-E10}
   (\omega_u)^n
  = \psi^s (z, u) \omega^n \;\; \mbox{in $M$}
\end{equation}
where $\psi^s (z, u) = (1-s) e^u + s \psi (z, u)$. Set
\[ S := \{s \in [0,1]: \mbox{equation~\eqref{gblq-E10} is solvable in
           $\cH \cap C^{2, \alpha} (M)$}\} \]
and let $u^s \in \cH \cap C^{2, \alpha} (M)$ be the unique solution of
\eqref{gblq-E10} for $s \in S$. Obviously $S \neq \emptyset$ as
$0 \in S$ with $u^0 = 0$. Moreover, by the $C^{2, \alpha}$ estimates
we see that $S$ is closed. We need to show that $S$ is also open in
and therefore equal to $[0, 1]$; $u^1$ is then the desired solution.

Let $s \in S$ and let $\Delta^s$ denote the Laplace operator of
$(M, \omega_{u^s})$. In local coordinates,
\[ \Delta^s v = \fg^{i\bj} v_{i\bj} = \fg^{i\bj}  \partial_i \bpartial_j v \]
where $\{{\fg}^{i\bj}\} = \{\fg^{s}_{i\bj}\}^{-1}$ and
$\fg^{s}_{i\bj} = g_{i\bj} + u^s_{i\bj}$.
Note that $\Delta^s - \psi^s_u$,
where $\psi^{s}_u = \psi^{s}_u (\cdot, u^{s})$,
is the linearized operator of equation~\eqref{gblq-E10} at $u^s$, .
We wish to prove that  for any $\phi \in C^{\alpha} (M, \omega_{u^s})$
there exists a unique solution $v \in C^{2, \alpha} (M, \omega_{u^s})$
to the equation
\begin{equation}
\label{gblq-E30}
\Delta^s v - \psi^s_u v = \phi,
\end{equation}
which implies by the implicit function theorem that
$S$ contains a neighborhood of $s$ and hence is open in $[0,1]$, completing
the proof.

The proof follows a standard approach, using the Lax-Milgram theorem and
the Fredholm alternative. For completeness we include it here.

Let $\gamma > 0$ and define a bilinear form on the Sobolev space $H^1 (M, \omega_{u^s})$ by
\begin{equation}
\label{gblq-E40}
\begin{aligned}
B [v, w] := \,& \int_M [\langle \nabla v + v \mbox{tr} \tilde{T},
                                  \nabla w \rangle_{\omega_{u^s}}
                  + (\gamma + \psi^{s}_u) v w] (\omega_{u^{s}})^n \\
          = \,& \int_M [{\fg}^{i\bj} (v_i + v \tilde{T}_{ik}^k)w_{\bj}
               + (\gamma + \psi^{s}_u) v w] (\omega_{u^{s}})^n
\end{aligned}
\end{equation}
 where $\tilde{T}$ denotes the torsion of $\omega_{u^s}$ and
$\mbox{tr} \tilde{T}$ its trace. In local coordinates,
\[ \mbox{tr} \tilde{T} = \tilde{T}_{ik}^k d z_i
      = {\fg}^{k\bj} (g_{i\bj k} - g_{k\bj i}) dz_i \]
so it only depends on the second derivatives of $u$.

It is clear that for $\gamma > 0$ sufficiently large $B$ satisfies the
Lax-Milgram hypotheses, i.e,
\begin{equation}
\label{gblq-E50}
|B [v, w]| \leq C \|v\|_{H^1 (\omega^s)} \|w\|_{H^1 (\omega^s)}
\end{equation}
by the Schwarz inequality, and
\begin{equation}
\label{gblq-E60}
 B [v, v] \geq c_0 \|v\|_{H^1 (\omega^s)}^2,
   \;\; \forall \; v \in H^1 (M, \omega_{u^s})
\end{equation}
where $c_0$ is a positive constant independent of $s \in [0,1]$ since
$\psi_u > 0$, $|u^s|_{C^2 (M)} \leq C$ and $M$ is compact.
By the Lax-Milgram theorem, for any $\phi \in L^2 (M, \omega_{u^s})$
there is a unique $v \in H^1 (M, \omega_{u^s})$ satisfying
\begin{equation}
\label{gblq-E70}
B  [v, w] = \int_M  \phi w (\omega_{u^{s}})^n \;\;
    \forall \; w \in H^1 (M, \omega_{u^s}).
\end{equation}
On the other hand,
\begin{equation}
\label{gblq-E80}
    B [v, w] = \int_M (- \Delta^s v + \psi^s_u v + \gamma v) w
               (\omega_{u^{s}})^n
\end{equation}
by integration by parts. Thus $v$ is a weak solution to the equation
\begin{equation}
\label{gblq-E90}
 L_{\gamma} v := \Delta^s v - \psi^s_u v - \gamma v = \phi.
\end{equation}
We write $v = L_{\gamma}^{-1} \phi$.

By the Sobolev embedding theorem the linear operator
\[ K := \gamma L_{\gamma}^{-1}: L^2 (M, \omega_{u^s}) \rightarrow L^2 (M, \omega_{u^s}) \]
is compact. Note also that $v \in H^1 (M, \omega_{u^s})$ is a weak solution of
equation~\eqref{gblq-E30} if and only if
\begin{equation}
\label{gblq-E100}
 v - K v = \zeta
\end{equation}
where $\zeta = L_{\gamma}^{-1} \phi$. Indeed, \eqref{gblq-E30} is
equivalent to
\begin{equation}
\label{gblq-E105}
v = L_{\gamma}^{-1} (\gamma v + \phi) = \gamma L_{\gamma}^{-1} v
   + L_{\gamma}^{-1} \phi.
\end{equation}
Since the solution of equation~\eqref{gblq-E30}, if exists, is unique,
by the Fredholm alternative equation~\eqref{gblq-E100} is uniquely solvable
for any $\zeta \in L^2 (M, \omega_{u^s})$.
Consequently, for any $\phi \in L^2 (M, \omega_{u^s})$ there exists a unique
solution $v \in H^1 (M, \omega_{u^s})$ to equation~\eqref{gblq-E30}.
By the regularity theory of linear elliptic equations,  $v \in C^{2, \alpha} (M, \omega_{u^s})$
if  $\phi \in C^{\alpha} (M, \omega_{u^s})$.
This completes the proof.
\end{proof}

\subsection{The Dirichlet problem}
We now turn to the proof of Theorem~\ref{gblq-th20}. Let
\[ \cA_{\ul{u}} = \{v \in \cH: \mbox{$v \geq \ul{u}$ in $M$,
                    $v = \ul{u}$ on $\partial M$}\}. \]
By the maximum principle, $v \leq h$ on $\bM$ for all $v \in \cA_{\ul{u}}$
where $h$ satisfies $\Delta h + n = 0$ in $M$ and $h = \ul{u}$ on
$\partial M$. Therefore we have $C^0$ bounds for solutions
of the Dirichlet problem~\eqref{gblq-I10}-\eqref{gblq-I20} in $\cA_{\ul{u}}$.
The proof of existence of such solutions then follows that of Theorem~1.1 in
\cite{Guan98a}; so is omitted here.

\begin{proof}[Proof of Theorem~\ref{gblq-th50}]
As we only assume $\psi \geq 0$, equation~\eqref{gblq-I100} is
degenerate. So we need to approximate it by nondegenerate equations.
Since $\omega_{\phi} > 0$ and $M$ is compact, there is $\varepsilon_0 > 0$
such that $\omega_{\phi} \geq \varepsilon_0 \omega$, and therefore
$(\omega_{\phi})^n \geq \varepsilon_0^n \omega^{n}$.

For $\varepsilon \in (0, \varepsilon_0]$ let $\psi^{\varepsilon}$ be a
smooth function such that
\[ \sup \Big\{\psi - \varepsilon, \frac{\varepsilon^n}{2}\Big\}
\leq \psi^{\varepsilon} \leq \sup \{\psi, \varepsilon^n \} \]
and consider the approximating problem
 \begin{equation}
\label{cma-K700'}
 \left\{ \begin{aligned}
 & (\omega_u)^n = \psi^\varepsilon \omega^{n} \;\; \mbox{in $\bM$}, \\
 & u = \phi \;\;  \mbox{on $\partial M$}.
  \end{aligned} \right.
\end{equation}
Note that $\phi$ is a subsolution of \eqref{cma-K700'} when
$0 < \varepsilon \leq \varepsilon_0$.
By Theorem~\ref{gblq-th20} there is a unique
solution $u^{\varepsilon} \in C^{2,\alpha} (\bM)$ of
\eqref{cma-K700'} with $u^{\varepsilon} \geq \phi$ on $\bM$
for $\varepsilon \in (0, \varepsilon_0]$.

By Theorem~\ref{gblq-th10} it is easy to see that
 \begin{equation}
\label{cma-K730}
 |u^{\varepsilon}|_{C^1(\bM)} \leq C_1, \; \sup_{M} \Delta u^{\varepsilon} \leq
C_2 (1+ \sup_{\partial M} \Delta u^{\varepsilon}),
\; \mbox{independent of $\varepsilon$}.
\end{equation}

On the boundary $\partial M$, the estimates in Section~\ref{gblq-B} for the pure
tangential and mixed tangential-normal second derivatives are independent of
$\varepsilon$, i.e.,
\begin{equation}
\label{cma-K750}
 |u^{\varepsilon}_{\xi \eta}|, \;  |u^{\varepsilon}_{\xi \nu}| \leq C_3,
\; \forall \; \xi, \eta  \in T \partial M, |\xi|, |\eta| = 1
\; \mbox{independent of $\varepsilon$}.
\end{equation}
where $\nu$ is the unit normal to $\partial M$.
For the estimate of the double normal derivative
$u^{\varepsilon}_{\nu \nu}$, note that $\partial M = N \times \partial S$
and $T_{C} \partial M = T N$; this is the only place we need the assumption
 $M = N \times S$ so Theorem~\ref{gblq-th50} actually holds for local
 product spaces.
 So
\begin{equation}
\label{cma-K740}
1 + u^{\varepsilon}_{\xi \bar{\xi}} = 1 + \phi_{\xi \bar{\xi}} \geq c_0 \;\;
\forall \; \xi \in T_{C} \partial M = T N, \; |\xi| = 1.
\end{equation}
where $c_0$ depends only on $\phi$.
From the proof in Section~\ref{gblq-B} we see that
\begin{equation}
\label{cma-K760}
 |u^{\varepsilon}_{\nu \nu}| \leq C,
\;\; \mbox{independent of $\varepsilon$ on $\partial M$}.
\end{equation}

Finally, from $\sup_M |\Delta u^{\varepsilon}| \leq C$ we see that
$|u^{\varepsilon}|_{C^{1, \alpha} (\bM)}$ is bounded for any
$\alpha \in (0, 1)$. Taking a convergent subsequence we obtain a
solution $u \in C^{1, \alpha} (\bM)$
of \eqref{gblq-I100} with the desired properties.
By Remark~\ref{gblq-remark-R10}, $u \in C^{1, 1} (\bM)$ when
$M$ has nonnegative bisectional curvature.
\end{proof}

\bigskip

\section{Geodesics in the space of Hermitian metrics.}
\label{gblq-S}
\setcounter{equation}{0}
\medskip

Let $(M, g)$ be a compact Hermitian manifold without boundary. The space of
Hermitian metrics
\begin{equation}
\label{cma-K220'}
 \cH = \{\phi \in C^{2} (M): \omega_{\phi} > 0\}
\end{equation}
is an open subset of $C^2 (M)$. The tangent space
$T_{\phi} \cH$ of $\cH$ at $\phi \in \cH$ is naturally identified to
$C^2 (M)$. Following Mabuchi~\cite{Mabuchi87}, Semmes~\cite{Semmes92}
and Donaldson~\cite{Donaldson99} who considered the K\"ahler case, we define
\begin{equation}
\label{cma-K620}
 \langle\xi, \eta\rangle_{\phi} = \int_M \xi \eta \, (\omega_{\phi})^n,
 \;\; \xi, \eta \in T_{\phi} \cH.
\end{equation}
Accordingly, the length of a regular curve
$\varphi: [0, 1] \rightarrow \cH$ is defined to be
 \begin{equation}
\label{cma-K630}
 L (\varphi) =
 \int_0^1 \langle\dot{\varphi}, \dot{\varphi}\rangle_{\varphi}^{\frac{1}{2}} dt.
 \end{equation}
Henceforth $\dot{\varphi} = \partial \varphi/\partial t$
and $\ddot{\varphi} = \partial^2 \varphi/\partial t^2$.
The geodesic equation takes the form
 \begin{equation}
\label{cma-K640}
 \ddot{\varphi} - |\nabla \dot{\varphi}|_{\varphi}^2 = 0,
 \end{equation}
 or in local coordinates
 \begin{equation}
\label{cma-K640'}
 \ddot{\varphi} - g (\varphi)^{j\bk} \dot{\varphi}_{z_j} \dot{\varphi}_{\bz_k}
 = 0.
 \end{equation}
Here $\{g (\varphi)^{j\bk}\}$ is the inverse matrix of
$\{g (\varphi)_{j\bk}\} = \{g_{j\bk} + \varphi_{j\bk}\}$.

It was observed by Donaldson~\cite{Donaldson99}, Mabuchi~\cite{Mabuchi87}
and Semmes~\cite{Semmes92} that the geodesic equation~\eqref{cma-K640}
reduces to a homogeneous complex Monge-Amp\`ere equation in $M \times A$
where $A = [0,1] \times \bfS^1$.
Let
\[ w = z_{n+1} = t + \sqrt{-1} s \]
 be a local coordinate of $A$.
We may view a smooth curve $\varphi$ in $\cH$ as a function on
$M \times [0, 1]$ and therefore a rotation-invariant function
(constant in $s$) on $M \times A$. Clearly,
\[ \dot{\varphi} = \frac{\partial \varphi}{\partial t}
   = 2 \frac{\partial \varphi}{\partial w}
   = 2 \frac{\partial \varphi}{\partial \bw}, \;\;
\ddot{\varphi} = \frac{\partial^2 \varphi}{\partial t^2}
   = 4 \frac{\partial^2 \varphi}{\partial w \partial \bw}. \]
Therefore,
 \begin{equation}
\label{cma-K660}
 \begin{aligned}
  \det \,& \begin{bmatrix}
  &                         &          &\varphi_{1\bw} \\
  & (g (\varphi)_{j\bar{k}})&                   & \vdots \\
  &                         &          &\varphi_{n\bar{w}} \\
\varphi_{w\bar{1}} & \cdots & \varphi_{w\bar{n}} & \varphi_{w\bw}
\end{bmatrix}  \\
  = \,& \frac{1}{4} \det (g (\varphi)_{i\bar{j}}) \cdot
 \det \begin{bmatrix}
  &                         &          & g (\varphi)^{k\bar{1}} \dot{\varphi}_{k} \\
  &  I                      &          & \vdots \\
  &                         &          & g (\varphi)^{k\bar{n}} \dot{\varphi}_{k} \\
\dot{\varphi}_{\bar{1}} & \cdots & \dot{\varphi}_{\bar{n}} &
\ddot{\varphi}
\end{bmatrix} \\
= &  \frac{1}{4} \det (g (\varphi)_{i\bar{j}}) \cdot
   (\ddot{\varphi} - g (\varphi)^{j\bk} \dot{\varphi}_{z_j}
   \dot{\varphi}_{\bz_k}).
 \end{aligned}
\end{equation}
So a geodesic $\varphi$ in $\cH$ satisfies
 \begin{equation}
\label{cma-K670}
 (\tilde{\omega}_{\varphi})^{n+1} \equiv \Big(\tilde{\omega} + \frac{\sqrt{-1}}{2} \partial
\bar{\partial} \varphi\Big)^{n+1} = 0 \;\;\; \mbox{in $M \times A$}
\end{equation}
where
 \begin{equation}
\label{cma-K680}
 \tilde{\omega} = \omega + \frac{\sqrt{-1}}{2} \partial \bar{\partial} |w|^2
  = \frac{\sqrt{-1}}{2} \Big(\sum_{j, k \leq n}  g_{j\bk} dz_j \wedge
  d\bz_k + dw \wedge d\bw\Big)
\end{equation}
is the lift of $\omega$ to $M \times A$.

Conversely, if $\varphi \in C^2 (M \times A)$ is a
rotation-invariant solution of \eqref{cma-K670} such that
\begin{equation}
\label{cma-K690}
 \varphi (\cdot, w) \in \cH, \;\; \forall \; w \in A,
 \end{equation}
then $\varphi$ is a geodesic in $\cH$.

In the K\"ahler case, Donaldson~\cite{Donaldson99} conjectured that
$\cH^{\infty} \equiv \cH \cap C^{\infty} (M)$ is geodesically convex,
i.e., any two functions in $\cH^{\infty}$ can be connected by a smooth
geodesic.
More precisely,

\begin{conjecture}[Donaldson~\cite{Donaldson99}]
\label{Donaldson-C99o}
Let $M$ be a compact K\"ahler manifold without boundary and
$\rho \in C^{\infty} (M \times \partial A)$ such that
$\rho (\cdot, w) \in \cH$ for $w \in \partial A$. Then there exists
a unique solution $\varphi$  of the Monge-Amp\`ere equation
\eqref{cma-K670} satisfying \eqref{cma-K690} and the boundary
condition $\varphi = \rho$.
\end{conjecture}

The uniqueness was proved by Donaldson~\cite{Donaldson99} as a
consequence of the maximum principle. In \cite{Chen00}, X.-X. Chen
obtained the existence of a weak solution
with $\Delta \varphi \in L^{\infty} (M \times A)$; see also the recent
work of Blocki~\cite{Blocki} who proved that the solution is in
$C^{1,1} (M \times A)$ when $M$ has nonnegative bisectional curvature.
As a corollary of Theorem~\ref{gblq-th50} these results can be extended
to the Hermitian case.

\begin{theorem}
\label{cma-K-thm30}
Let $M$ be a compact Hermitian manifold without boundary.
and let
$\varphi_0, \varphi_1 \in \cH \cap C^4 (M)$.
 There exists a unique (weak) solution $\varphi \in C^{1,\alpha} (M \times
 A)$, $\forall \; 0 < \alpha < 1$, with $\tilde{\omega}_{\varphi} \geq 0$
 and $\Delta \varphi \in L^{\infty} (M \times A)$
 of the Dirichlet problem
  \begin{equation}
\label{cma-K700}
 \left\{ \begin{aligned}
 & (\tilde{\omega}_{\varphi})^{n+1}  = 0 \;\; \mbox{in $M \times A$} \\
 & \varphi = \varphi_0 \;\; \mbox{on $M \times \Gamma_0$}, \\
 &  \varphi = \varphi_1 \;\; \mbox{on $M \times \Gamma_1$}
  \end{aligned} \right.
\end{equation}
where $\Gamma_0 = \partial A|_{t=0}$, $\Gamma_1 = \partial A|_{t=1}$.
Moreover, $\varphi \in C^{1,1} (M \times A)$ if $M$ has nonnegative
bisectional curvature.
\end{theorem}

\begin{proof}
In order to apply Theorem~\ref{gblq-th50} to the Dirichlet
problem~\eqref{cma-K700} we only need to construct a strict subsolution.
This is easily done for the
annulus $A = [0,1] \times \bfS^1$. Let
\[ \ul{\varphi} = (1 - t) \varphi_0 + t \varphi_1 + K (t^2 - t). \]
Since $\varphi_0, \varphi_1 \in \cH (\omega)$ we see that
$\tilde{\omega}_{\ul{\varphi}} > 0$ and
$(\tilde{\omega}_{\ul{\varphi}})^{n+1} \geq 1$
for $K$ sufficiently large.
\end{proof}

\begin{remark}
By the uniqueness $\varphi$ is rotation invariant (i.e., independent of $s$).
\end{remark}


\bigskip

\small
\begin{thebibliography}{99}

\bibitem{Aubin78}
T. Aubin,
{\em \'Equations du type Monge-Amp\`ere sur les vari\'et\'es k\"ahl\'eriennes
compactes}, (French)  Bull. Sci. Math. (2) {\bf 102} (1978), 63--95.

\bibitem{BF79}
E. Bedford and J. E. Fornaess,
{\em Counterexamples to regularity for the complex Monge-Amp\`ere equation},
Invent. Math. {\bf 50} (1979), 129--134.

\bibitem{BT76}
E. Bedford and B. A. Taylor,
{\em The Dirichlet problem for a complex Monge-Amp\`ere equation},
Invent. Math. {\bf 37} (1976), 1--44.

\bibitem{BT79}
E. Bedford and B. A. Taylor,
{\em  Variational properties of the complex Monge-Amp\`ere equation,
II. Intrinsic norms},
Amer. J. Math., {\bf 101} (1979), 1131--1166.

\bibitem{BT82}
E. Bedford and B. A. Taylor, {\em A new capacity for
plurisubharmonic functions}, Acta Math. {\bf 149} (1982), 1--40.

\bibitem{Blocki05}
Z. Blocki,
{\em On uniform estimate in Calabi-Yau theorem}, Proc. SCV2004, Beijing,
Science in China Series A Mathematics 48 Supp. (2005), 244--247.

\bibitem{Blocki09}
Z. Blocki,
{\em A gradient estimate in the Calabi-Yau theorem},
Math. Annalen {\bf 344} (2009), 317--327.

\bibitem{Blocki}
Z. Blocki,
{\em On geodesics in the space of K\"ahler metrics},
preprint.

\bibitem{CKNS}
L. A. Caffarelli, J. J. Kohn, L. Nirenberg and J. Spruck, {\em The
Dirichlet problem for nonlinear second-order elliptic equations II.
Complex Monge-Amp\`{e}re and uniformly elliptic equations},
{ Comm. Pure Applied Math.} {\bf 38} (1985), 209--252.

\bibitem{CNS1}
L. A. Caffarelli, L. Nirenberg and J. Spruck,
{\em The Dirichlet problem for nonlinear second-order elliptic equations I.
Monge-Amp\`{e}re equations},
{ Comm. Pure Applied Math.} {\bf 37} (1984), 369--402.

\bibitem{Calabi56}
E. Calabi,
{\em The space of K\"aler metrics},
Proc. ICM, Amsterdam 1954, Vol. 2, 206--207, North-Holland, Amsterdam, 1956.

\bibitem{CP92}
U. Cegrell and L. Persson,
{\em The Dirichlet problem for the complex Monge-Amp\"ere operator: stability
in $L^2$},  Michigan Math. J. {\bf 39} (1992), 145--151.

\bibitem{Chen00}
X.-X. Chen,
{\em The space of K\"ahler metrics},
J. Differential Geom. {\bf 56} (2000), 189--234.

\bibitem{CW}
Y.-Z. Chen and L.-C. Wu,
{\em Second Order Elliptic Equations and Elliptic Systems}, Amec. Math. Soc.,
Providence, RI, 1998.


\bibitem{CLN69}
S. S. Chern, H. I. Levine, L. Nirenberg,
{\em Intrinsic norms on a complex manifold},
1969 Global Analysis (Papers in Honor of K. Kodaira)
pp. 119--139, Univ. Tokyo Press, Tokyo.

\bibitem{Donaldson99}
 S. K. Donaldson,
{\em Symmetric spaces, K\"ahler geometry and Hamiltonian dynamics},
Northern California Symplectic Geometry Seminar, 13--33,
Amer. Math. Soc. Transl. Ser. 2, {\bf 196},
Amer. Math. Soc., Providence, RI, 1999.

\bibitem{Evans}
L. C. Evans,
{\em Classical solutions of fully nonlinear, convex, second order elliptic
equations},
{ Comm. Pure Applied Math.} {\bf 35} (1982), 333--363.

\bibitem{GS80}
T. W Gamelin and N. Sibony,
{\em Subharmonicity for uniform algebras},
J. Funct. Anal. {\bf 35} (1980), 64--108.

\bibitem{Guan98a}
B. Guan,
{\em The Dirichlet problem for Monge-Amp\`ere equations in non-convex
domains and spacelike hypersurfaces of constant Gauss curvature},
Trans. Amer. Math. Soc. {\bf 350} (1998), 4955--4971.

\bibitem{Guan98b}
B. Guan, {\em The Dirichlet problem for complex Monge-Amp\`ere
equations and regularity of the pluri-complex Green function},
Comm. Anal. Geom. {\bf 6} (1998), 687--703. {\em A correction},
 {\bf 8} (2000), 213--218.

\bibitem{GL2}
B. Guan and Q. Li,
{\em Totally real submanifolds and the homogeneous complex Monge-Amp\`ere
equation}, in preparation.

\bibitem{GL96}
B. Guan and Y.-Y. Li,
{\em Monge-Amp\`ere equations on Riemannian manifolds},
J. Differential Equations {\bf 132} (1996), 126--139.

\bibitem{GS93}
B.~Guan and J.~Spruck,
{\em Boundary value problem on $\bfS^n$ for surfaces of constant Gauss
 curvature},
{ Annals of Math.} {\bf 138} (1993), 601--624.

\bibitem{GuanPF02}
P.-F. Guan,
{\em Extremal  functions related to intrinsic norms},
Annals of Math. {\bf 156} (2002), 197--211.

\bibitem{GuanPF08}
P.-F. Guan,
{\em Remarks on the homogeneous complex Monge-Amp\`ere equation related to
the Chern-Levine-Nirenberg conjecture}, preprint.

\bibitem{GuanPF}
P.-F. Guan,
{\em A gradient estimate for complex Monge-Amp\`ere equation},
unpublished.

\bibitem{Guillemin-Stenzel91}
V. Guillemin and M. Stenzel,
{\em Grauert tubes and the homogeneous Monge-Amp\`ere equation},
J. Differential Geom. \textbf{34} (1991), 561--570.

\bibitem{Guillemin-Stenzel92}
V. Guillemin and M. Stenzel,
{\em Grauert tubes and the homogeneous Monge-Amp\`ere equation II},
J. Differential Geom. \textbf{35} (1992), 627--641.

\bibitem{HRS}
D. Hoffman, H. Rosenberg and J. Spruck,
{\em Boundary value problems for surfaces of constant Gauss curvature},
{ Comm. Pure Applied Math.} {\bf 45} (1992), 1051--1062.

\bibitem{Kazdan78}
J. L. Kazdan,
{\em A remark on the preceding paper of S. T. Yau: ``On the Ricci curvature
of a compact K\"ahler manifold and the complex Monge-Amp\`ere equation. I''},
Comm. Pure Appl. Math. {\bf 31} (1978), 413--414.

\bibitem{Kolodziej98}
S. Kolodziej,
{\em The complex Monge-Amp\`ere equation}, Acta Math. {\bf 180} (1998), 69-117.

\bibitem{Krylov82}
N. V.~Krylov,
{\em Boundedly inhomogeneous elliptic and parabolic equations},
(Russian)  Izv. Akad. Nauk SSSR Ser. Mat. {\bf 46} (1982), 487--523, 670.
English translation: Math. USSR Izv. {\bf 22} (1984), 67--98.

\bibitem{Krylov83}
N. V.~Krylov,
{\em Boundedly inhomogeneous elliptic and parabolic equations in a domain}
(Russian)  Izv. Akad. Nauk SSSR Ser. Mat. {\bf 47} (1983), 75--108.

\bibitem{Krylov}
N. V.~Krylov,
{\em Smoothness of the payoff function for a controllable diffusion process
in a domain} (Russian)  Izv. Akad. Nauk SSSR Ser. Mat.  {\bf 53}  (1989),
66-96;  translation in  Math. USSR-Izv. {\bf 34} (1990), 65--95.

\bibitem{Lempert-Szoke}
L. Lempert and R. Sz\"oke
{\em Global solutions of the homogeneous complex Monge-Amp\`ere equation
and complex structures on the tangent bundle of Riemannian manifolds},
Math. Ann. {\bf 290} (1991), 689--712.

\bibitem{Mabuchi87}
T. Mabuchi,
{\em Some symplectic geometry on compact K\"ahler manifolds. I},
 Osaka J. Math. {\bf 24} (1987), 227--252.

\bibitem{Patrizio-Wong}
G. Patrizio and P.-M. Wong,
{\em Stein manifolds with compact symmetric center},
Math. Ann. {\bf 289} (1991), 355--382.

\bibitem{Semmes92}
S. Semmes,
{\em Complex Monge-Amp\`ere and symplectic manifolds},
Amer. J. Math.  {\bf 114} (1992), 495--550.

\bibitem{Siu87}
Y.-T. Siu,
{\em Lectures on Hermitian-Einstein metrics for stable bundles and
K\"ahler-Einstein metrics},
DMV Seminar, 8. Birkh\"auser Verlag, Basel, 1987. 171 pp.

\bibitem{ST}
J. Street and G. Tian,
{\em Hermitian curvature flow}, preprint.

\bibitem{Tian00}
G. Tian,
{\em Canonical metrics in K\"ahler geometry. Notes taken by Meike Akveld.}
Lectures in Mathematics ETH Z\"urich. Birkh\"auser Verlag, Basel, 2000.
vi+101 pp.

\bibitem{Wong}
P.-M. Wong,
{\em Geometry of the complex homogeneous Monge-Amp\`ere equation},
Invent. Math. \textbf{67} (1982), 349--373.

\bibitem{Yau78}
S.-T. Yau,
{\em On the Ricci curvature of a compact K\"ahler manifold and the complex
Monge-Amp\`ere equation. I.}
Comm. Pure Appl. Math. {\bf 31} (1978), no. 3, 339--411.

\bibitem{Zhang}
X.-W. Zhang,
{\em  A gradient estimate for complex Momnge-Amp\`ere equation},
preprint.


\end{thebibliography}
\end{document}